\documentclass[11pt]{article}

\usepackage{amssymb}
\usepackage{amsmath,amsfonts,mathtools}
\usepackage{amsthm}
\usepackage{latexsym}
\usepackage{mathrsfs}
\usepackage{color}
\usepackage{float}
\usepackage{enumitem}
\usepackage{algorithm, algorithmic}
\allowdisplaybreaks[3]

\usepackage{bbm}
\usepackage{microtype}
\usepackage[a4paper]{geometry}
\usepackage{pifont}
\usepackage{booktabs}

\usepackage{xcolor}

\newtheorem{lemma}{Lemma}
\newtheorem{theorem}{Theorem}
\newtheorem{corollary}{Corollary}
\newtheorem{definition}{Definition}
\newtheorem{remark}{Remark}
\newtheorem{assumption}{Assumption}

\newcommand{\beq}{\begin{equation}}
\newcommand{\eeq}{\end{equation}}
\newcommand{\nn}{\nonumber}

\def\bc{{\bar c}}
\def\bF{{\bar F}}
\def\bh{{\bar h}}
\def\bi{{\bar i}}
\def\bk{{\bar k}}
\def\bu{{\bar u}}

\def\by{{\bar y}}
\def\bz{{\bar z}}
\def\bPsi{{\bar \Psi}}
\def\bR{{\mathbb R}}
\def\cB{{\mathcal B}}

\def\cL{{\mathcal L}}
\def\cX{{\mathcal X}}
\def\cY{{\mathcal Y}}
\def\cS{{\mathcal S}}
\def\cN{{\mathcal N}}
\def\cU{{\mathcal U}}

\def\dist{{\rm dist}}
\def\prox{{\rm prox}}
\def\dom{{\rm dom}}

\def\rc{{\rm c}}
\def\tz{{\tilde z}}
\def\rC{{\rm C}}

\title{A first-order method for constrained nonconvex-nonconcave minimax optimization}
\author{
Zhaosong Lu
\thanks{
Department of Industrial and Systems Engineering, University of Minnesota, USA (email: {\tt zhaosong@umn.edu}, {\tt wan02269@umn.edu}). This work was partially supported by the Air Force Office of Scientific Research under Award FA9550-24-1-0343, the Office of Naval Research under Award N00014-24-1-2702, and the National Science Foundation under Awards 2211491 and 2435911.
}
\and
Xiangyuan Wang
\footnotemark[1]
}
\date{October 11, 2025 (Revised: May 21, 2026)}

\begin{document}

\maketitle

\begin{abstract}
We study a class of constrained nonconvex-nonconcave minimax optimization problems in which the inner maximization involves potentially complex constraints. Under the assumption that the inner problem of a novel lifted minimax reformulation satisfies a local Kurdyka-\L{}ojasiewicz (KL) condition, we show that the maximal function of the original problem enjoys a local {generalized} H\"{o}lder smoothness property. We also propose a sequential convex programming (SCP) method for solving constrained optimization problems and establish its convergence rate under a local KL condition. Leveraging these results, we develop an inexact proximal gradient method for the original minimax problem, where the inexact gradient of the maximal function is computed via the SCP method applied to a locally KL-structured subproblem. Finally, we establish complexity guarantees for the proposed method in computing an approximate stationary point of the original minimax problem.
\end{abstract}

\noindent{\bf Keywords:} constrained nonconvex-nonconcave minimax optimization, local KL condition, local {generalized} H\"{o}lder smoothness, sequential convex programming method, inexact proximal gradient method, first-order oracle complexity

\medskip

\noindent {\bf Mathematics Subject Classification:} 90C26, 90C30, 90C47, 90C99, 65K05 

\section{Introduction}\label{sec:intro}

In this paper, we consider a class of constrained nonconvex-nonconcave minimax optimization problems of the form
\beq \label{intro:problem}
\min_x \max_{c(y) \leq 0} \left\{ f(x, y) + p(x) - q(y) \right\},
\eeq
where $f$ is a smooth function that is possibly nonconvex in $x$ and nonconcave in $y$, $p$ and $q$ are simple closed convex functions, and $c$ is a smooth mapping. This problem arises in many applications and also appears as a subproblem when solving more general constrained minimax problems of the form
 \beq \label{general-problem}
\min_{d(x) \leq 0} \max_{c(y) \leq 0} \left\{ f(x, y) + p(x) - q(y) \right\},
\eeq
where $d$ is a smooth mapping. In fact, by applying a penalty approach, one can naturally convert \eqref{general-problem} into a sequence of subproblems
\[
\min_{x} \max_{c(y) \leq 0} \left\{ f(x, y)+\rho_k \|[d(x)]_+\|^2 + p(x) - q(y) \right\},
\]
which are clearly in the form of \eqref{intro:problem}, where $0<\rho_k \to \infty$ and $u_+=\max\{u,0\}$ for any vector $u$.

In recent years, considerable attention has been devoted to unconstrained nonconvex-nonconcave minimax problems of the following form:
\beq \label{intro_prSimpler}
\min_x \max_y \left\{ f(x, y) + p(x) - q(y) \right\}.
\eeq
This class of problems arises in a wide range of applications in machine learning and operations research, including generative adversarial networks~\cite{arjovsky2017wasserstein, goodfellow2020generative}, reinforcement learning~\cite{dai2018sbeed, omidshafiei2017deep}, adversarial training~\cite{madry2017towards, sinha2017certifying}, and distributionally robust optimization~\cite{bertsimas2011theory, blanchet2025distributionally, rahimian2022frameworks}. Significant progress has been made in solving~\eqref{intro_prSimpler} under additional structural assumptions. For instance, several works study the special case with $q=0$ and assume that the inner maximization problem in~\eqref{intro_prSimpler} satisfies a global Polyak–\L{}ojasiewicz (PL) condition, which is generally weaker than strong concavity. Under this assumption, gradient descent–ascent type methods have been developed, and complexity guarantees have been established for obtaining approximate stationary points (see, 
e.g.,  \cite{huang2023enhanced,nouiehed2019solving, xu2023zeroth,yang2022faster}). In addition, first-order methods have been proposed for problem~\eqref{intro_prSimpler} from the perspective of variational inequalities, typically assuming the existence of a weak Minty variational inequality solution (see, e.g.,~\cite{bohm2022solving,cai2022accelerated,liu2021first,pethick2023escaping}).

More recently, \cite{li2025nonsmooth, zheng2025doubly, zheng2023universal} studied problem~\eqref{intro_prSimpler} under a global KL condition, where $p$ and $q$ are indicator functions of simple convex compact sets. This setting generalizes that of~\cite{huang2023enhanced, nouiehed2019solving, xu2023zeroth, yang2022faster}, since the KL condition extends the PL condition (the latter corresponding to the KL condition with exponent $1/2$).  However, requiring the KL property to hold globally is often too restrictive in practice. To address this, \cite{lu2025first} considered problem~\eqref{intro_prSimpler} with $p$ and $q$ being simple closed convex functions under a local KL condition. Specifically, for each fixed outer variable $x \in \dom\,p$, the KL property is assumed to hold only on a level set of the inner variable $y$, where this level set may depend on $x$ and may shrink as $x$ approaches a stationary point of problem~\eqref{intro_prSimpler}. Under this weaker assumption, a local {generalized} H\"{o}lder smoothness property of the associated maximal function was established. Leveraging this property, an inexact proximal gradient method was developed, in which the inexact gradient of the maximal function is computed by applying a proximal gradient method to a locally KL-structured subproblem. Complexity guarantees were then established for finding an approximate stationary point of problem~\eqref{intro_prSimpler}.

Despite recent advances, existing results primarily focus on nonconvex-nonconcave minimax problems with unconstrained inner maximization. To the best of our knowledge, no algorithmic framework has been developed for the constrained counterpart~\eqref{intro:problem}, where the inner maximization problem involves potentially complex constraints. In this paper, we study problem~\eqref{intro:problem} under the assumption that a novel \emph{lifted minimax} problem---equivalent to \eqref{intro:problem}---satisfies a local KL condition analogous to that considered for problem~\eqref{intro_prSimpler} (see Assumption~\ref{ass:Lip_Smo_KL}). We establish that the maximal function $F^*(x) := \max_{c(y) \leq 0} \{ f(x, y) - q(y) \}$ exhibits a local {generalized} H\"{o}lder smoothness property. In addition, for any fixed outer variable $x \in \dom\,p$, we propose a sequential convex programming (SCP) method to solve the inner maximization subproblem and establish its convergence rate under the local KL condition. Building on these results, we develop an inexact proximal gradient method for solving $\min_x \{ F^*(x) + p(x) \}$, which is equivalent to \eqref{intro:problem}. Specifically, given the current iterate $(x^k, y^{k-1})$, we apply the SCP method to approximately solve $\max_{c(y)\le 0} \{ f(x^k, y) - q(y)\}$ initialized at $y^{k-1}$, and obtain an approximate solution $y^k$. We then update $x^{k+1}$ via an inexact proximal gradient step, {using $-\nabla_x f(x^k, y^k)$ as the forward direction together with} a suitably chosen step size. Finally, we establish complexity guarantees for the proposed method in computing an approximate stationary point of problem~\eqref{intro:problem}.

The main contributions of this paper are summarized below.

\begin{itemize}
    \item We establish a local {generalized} H\"{o}lder smoothness property for the maximal function $F^*$ under a local KL condition imposed on a novel lifted minimax problem, which is crucial for developing a method for solving problem~\eqref{intro:problem}.

    \item We propose a sequential convex programming method  for solving a constrained optimization problem and establish its convergence rate under a local KL condition. This method serves as a subroutine for solving problem~\eqref{intro:problem}.

    \item We propose an inexact proximal gradient method for finding approximate stationary points of problem~\eqref{intro:problem}, and show that it achieves an \emph{iteration complexity} of $\widetilde{\mathcal{O}}\big(\epsilon^{-\max\{(1-\theta)^{-1},\, \theta^{-1}\sigma\}}\big)$, and a \emph{first-order oracle complexity} of $\widetilde{\mathcal{O}}\big(\epsilon^{-(1-\theta)^{-1}(2\theta^2 - 2\theta + 1)\max\{(1-\theta)^{-1},\, \theta^{-1}\sigma\}}\big)$, measured by the number of gradient evaluations, for finding an $\mathcal{O}(\epsilon)$-approximate stationary point of \eqref{intro:problem}, where $\theta$ and $\sigma$ are the parameters of the local KL condition.
\end{itemize}

The rest of this paper is organized as follows. Subsection~\ref{subsec:notation} introduces the notation, terminology, and assumptions used throughout the paper. In Section~\ref{sec:localKL_Holder}, we study the theoretical properties of problem~\eqref{intro:problem}. Section~\ref{sec:subsolver_SCP} presents a sequential convex programming method for a constrained optimization problem satisfying a local KL property. In Section~\ref{sec:FOD}, we propose an inexact proximal gradient method for solving problem~\eqref{intro:problem} and establish its complexity results. Section~\ref{sec:Numerical} presents  preliminary numerical results illustrating the performance of the proposed method. Finally, we provide the proof of the main results in Section \ref{sec:proof}.

\subsection{Notation, terminology, and assumptions} \label{subsec:notation}

The following notation will be used throughout the paper. Let $\mathbb{R}^n$ stand for the $n$-dimensional Euclidean space, and $\overline{\bR}=(-\infty, \infty]$. The standard inner product, $\ell_1$-norm, $\ell_\infty$-norm, and Euclidean norm are denoted by $\langle \cdot, \cdot \rangle$, $\|\cdot\|_1$, $\|\cdot\|_\infty$, and $\|\cdot\|$, respectively. For any two points $u, v \in \mathbb{R}^n$, the notation $[u, v]$ denotes the line segment connecting $u$ and $v$. Given a point $x$ and a closed set $S \subset \mathbb{R}^n$, let $\dist(x,S)$ stand for the distance from $x$ to $S$, and $\delta_{S}$ the indicator function of $S$.  The \emph{regular normal cone} and the \emph{normal cone} (i.e., the \emph{limiting normal cone}) of $S$ at $x \in S$ are denoted by $\widehat{\cN}_{S}(x)$ and $\cN_{S}(x)$, respectively (see \cite[Definition 6.3]{rockafellar2009variational}). The closed ball centered at $x \in \mathbb{R}^n$ with radius $r$ is denoted by $\cB(x, r)$. In addition, $\operatorname{conv}(\cdot)$, $\operatorname{aff}(\cdot)$, and $\operatorname{int}(\cdot)$ denote the convex hull, affine hull, and interior of the associated set, respectively.

A function $\phi: \cX \subseteq \mathbb{R}^n\to\mathbb{R}$ is called \emph{$L_\phi$-Lipschitz continuous} on $\cX$ if $\lvert \phi(x)-\phi(y)\rvert\le L_\phi\|x-y\|$ for all $x,y\in  \cX$, and \emph{$L_{\nabla \phi}$-smooth} on $ \cX$ if $\|\nabla \phi(x)-\nabla \phi(y)\|\le L_{\nabla \phi}\|x-y\|$ for all $x,y \in \cX$. For a closed convex function $p: \bR^n \to \overline{\bR}$, the \emph{proximal operator} associated with $p$ is defined as
\[
\prox_p(x) := \mathop{\arg\min}_{x' \in \bR^n} \Big\{ p(x') + \frac{1}{2} \|x - x'\|^2 \Big\}.
\]

For a function $\phi : \mathbb{R}^n \to \overline{\bR}$, its \emph{domain} is defined as $\dom\,\phi = \{x : \phi(x) < \infty\}$. Such $\phi$ is called \emph{proper} if $\dom\,\phi \neq \emptyset$, and it is called \emph{closed} or \emph{lower semicontinuous} if $\liminf_{z\to x} \phi(z)\geq \phi(x)$ holds for all $x\in\bR^n$. The \emph{regular subdifferential} (see, e.g., \cite[Definition 8.3(a)]{rockafellar2009variational}) of a proper closed function $\phi$ at $x \in \mathrm{dom}\, \phi$ is defined as
\[
\widehat{\partial} \phi(x) := \left\{ v \in \mathbb{R}^n : \liminf_{z \to x,\, z \neq x} \frac{\phi(z) - \phi(x) - \langle v, z - x \rangle}{\|z - x\|} \geq 0 \right\}.
\]
Let $z \xrightarrow{\phi} x$ denote $z \to x$ and $\phi(z) \to \phi(x)$. The \emph{limiting subdifferential} (see, e.g., \cite[Definition 8.3(b)]{rockafellar2009variational}) of a proper closed function $\phi$  at $x \in \mathrm{dom}\, \phi$ is defined as
\[
\partial \phi(x) := \left\{ v \in \mathbb{R}^n : \exists\, x^k \xrightarrow{\phi} x,\ v^k \to v\ \text{with } v^k \in \widehat{\partial} \phi(x^k) \right\}.
\]
We use $\partial_{x_i} \phi$ to denote the limiting subdifferential with respect to $x_i$. For an upper semicontinuous function $\phi$, its limiting subdifferential is defined as $\partial \phi = -\partial(-\phi)$. If $\phi$ is continuously differentiable, then $\partial \phi$ coincides with the gradient $\nabla \phi$. Besides, if $\phi$ is convex, then $\partial \phi$ corresponds to the classical convex subdifferential. It is well-known that $\partial (\phi_1 + \phi_2)(x) = \nabla \phi_1(x) + \partial \phi_2(x)$ if $\phi_1$ is continuously differentiable at $x$ and $\phi_2$ is lower or upper semicontinuous at $x$ (see, e.g., \cite[Exercise 8.8(c)]{rockafellar2009variational}).

{Suppose that $\phi$ is a locally Lipschitz continuous function on $\cX$. The \emph{restricted Clarke subdifferential} of $\phi$ with respect to $\mathcal{X}$, denoted by $\partial^\rC_\cX \phi(x)$,  is introduced in \cite{lu2025first} and defined as
\[
\partial^\rC_\cX \phi(x) := {\rm conv}\{v : \exists\, x^k \in \cX \rightarrow x  \text{ such that } \nabla \phi(x^k) \to v \} \quad \forall x\in \cX.
\]
When $\partial^\rC_\cX \phi(x)$ is a singleton, we denote its unique element by $\nabla^\rC_\cX \phi(x)$. Furthermore, suppose additionally that $\partial^\rC_\cX \phi(x)$ is a singleton for all $x \in \cX$. The function $\phi$ is said to be \emph{generalized H\"{o}lder smooth} on $\cX$ (see \cite[Definition 1]{lu2025first}) if there exist $L_1 \geq 0$, $L_2 \geq 0$, and  $\nu \in (0, 1]$ such that 
\[
\|\nabla^\rC_\cX \phi(x)-\nabla^\rC_\cX \phi(y)\|\le L_1 \|x-y\|+L_2\|x-y\|^\nu \quad \forall x,y\in  \cX.
\]}

For a closed function $\phi : \mathbb{R}^n \to \overline{\bR}$, the \emph{slope} of $\phi$ at $x\in\dom\,\phi$ is defined as
\beq \label{slope}
|\nabla \phi|(x):= \limsup_{z\to x} \frac{\big(\phi(x)-\phi(z)\big)_+}{\|x-z\|},
\eeq 
where $t_+=\max\{0,t\}$ for any $t\in\bR$. The \emph{limiting slope} of $\phi$ at $x\in\dom\,\phi$ is defined as 
\beq \label{limit-slope}
\overline{|\nabla \phi|}(x):= \liminf_{z \xrightarrow{\phi} x} |\nabla \phi|(z).
\eeq 
If $\phi$ is differentiable at $x$, $|\nabla \phi|(x)$ coincides with $\|\nabla \phi(x)\|$. When $\phi$ is a convex function, $|\nabla \phi|(x)$ reduces to $\dist(0,\partial \phi(x))$ for all $x\in\dom\,\phi$. For more details on the slope and limiting slope, see, for example,  \cite[Section 2]{drusvyatskiy2021nonsmooth}.

We now introduce the notion of an approximate stationary point for the problem $\min_x \phi(x)$, where $\phi$ is a closed function. Since the minimax problem~\eqref{intro:problem} can be viewed as a special case of this general problem, the following definition applies directly to~\eqref{intro:problem} as well.

\begin{definition}[{\bf $(r,\epsilon)$-stationary point}] \label{def-statpt}
Suppose $\phi$ is a closed function.  For any $\epsilon > 0$ and $r > 0$, a point $\bar{x}$ is called an $(r,\epsilon)$-stationary point of the problem $\min_x \phi(x)$ if $\bar{x}\in\dom\,\phi$ and $\dist(\bar{x}, \cX_\epsilon) \leq r$,
where $\cX_\epsilon = \{x \in \dom\,\phi: \dist(0, \partial \phi(x)) \le \epsilon\}$.
\end{definition}

{The above definition is closely related to the notion of an $\epsilon$-optimization-stationary point (see, e.g., \cite[Definition 3.1(a)]{li2025nonsmooth}). In particular, suppose that $\phi$ is a $\rho$-weakly convex function, i.e., $\phi(\cdot) + \rho\|\cdot\|^2/2$ is convex, and that $\bar{x}$ is an $\epsilon$-optimization-stationary point of the problem $\min_x \phi(x)$, that is, $\bar{x} \in \dom \phi$ and $\|\prox_{\phi/t}(\bar{x}) - \bar{x}\| \le \epsilon$ for some $t > \rho$. It can be verified that $\bar{x}$ is an $(\epsilon, t\epsilon)$-stationary point of $\min_x \phi(x)$. It should also be noted that when $\phi$ is locally Lipschitz continuous, any $(r,\epsilon)$-stationary point $\bar{x}$ of $\phi$ is also an $(r,\epsilon)$-Goldstein stationary point of $\phi$, that is, $\dist(0,\partial_r \phi(\bar{x})) \le \epsilon$, where
\[
\partial_r \phi(\bar{x}) := \operatorname{conv}\Big( \bigcup_{x \in \mathcal{B}(\bar{x}, r)} \partial \phi(x) \Big).
\]}

We next introduce additional notation for problem~\eqref{intro:problem}. For convenience, we define
\begin{align}
    &\cX := \dom\,p, \quad \cY := \dom\,q, \quad \cS := \cY \cap \{y : c(y) \leq 0\}, \label{intro_domain} \\
&  D_{\cY} = \max_{u, v \in \cY} \|u-v\|, \quad M_q = \max_{u, v \in \cY} \{q(u)-q(v)\},  \label{DY-Mq} \\
    &F(x, y) := f(x, y) - q(y) - \delta_{c(\cdot) \leq 0}(y), \quad F^*(x) =\max_{y} F(x, y), \quad Y^*(x) = \{y:F(x, y)=F^*(x)\}, \label{intro_F_Ystar} \\
&  \Psi(x) := F^*(x)+p(x), \quad \Psi^* := \min_x \Psi(x). \label{intro_Phi_star}
\end{align}
To study problem~\eqref{intro:problem}, we introduce the following \emph{lifted minimax} problem:
\beq \label{lift-minimax}
\min_x \max_{\bc(y,z) \leq 0} \left\{ f(x, y) + p(x) - q(y)- \delta_{\cY}(z) \right\},
\eeq
 where 
\beq \label{intro_cOver}
\bc=(\bc_1,\ldots,\bc_m) \ \ \text{with} \ \  \bc_j(y, z) :=  c_j(z) + \langle \nabla c_j(z), y - z \rangle + \frac{L_{c_j}}{2} \|y - z\|^2  \quad j=1,\ldots, m, 
\eeq
and $L_{c_j}$ is the Lipschitz smoothness constant of $c_j$ (see Assumption~\ref{ass:Lip_Smo_KL} below). Interestingly, the lifted minimax problem \eqref{lift-minimax} is equivalent to the original minimax problem  \eqref{intro:problem}  (see Lemma \ref{lem:equiv_Phi} and Remark \ref{minimax-equiv}). For notational convenience, we further define
\begin{align}
&\bF(x, y, z) := f(x, y) - q(y) - \delta_{\,\bc(\cdot, \cdot) \leq 0}(y, z) - \delta_{\cY}(z), \label{intro_potential} \\ 
& \bF^*(x) := \max_{y, z} \bF(x, y, z), \quad \overline{Y}^*(x) := \{(y, z): \bF(x, y, z) = \bF^*(x)\}, \quad \bPsi(x) := \bF^*(x) + p(x). \label{intro_PhiOver}
\end{align}
Consequently, problem~\eqref{lift-minimax} can be equivalently written as
\beq \label{equiv-prob}
\min_x \max_{y,z} \{\bF(x, y, z)+p(x)\}.
\eeq

In what follows, we introduce the assumptions for problem~\eqref{intro:problem} and its equivalent problem~\eqref{equiv-prob}. In particular, we assume that problem~\eqref{intro:problem} has at least one optimal solution and that the following assumption holds.

\begin{assumption} \label{ass:Lip_Smo_KL}
\begin{enumerate}[label=(\roman*)]
    \item For any fixed $y \in \cY$, the function $f(\cdot, y)$ is $L_f$-Lipschitz continuous on {$\cX$}. Moreover, the function $f: \bR^{n_1} \times \bR^{n_2} \to \bR$ is $L_{\nabla f}$-smooth on ${\cX} \times \cY$.
    
    \item The functions $p, q$ are proper closed convex, and the proximal operators of $p$ and $q$ can be computed exactly. Moreover, $\dom\,q$ (i.e., $\cY$) is compact. {In addition, we assume that $\operatorname{aff}(\cX) = \bR^{n_1}$.\footnote{{The assumption $\operatorname{aff}(\mathcal{X}) = \mathbb{R}^{n_1}$ is imposed merely for convenience. Without this assumption, the results of the paper remain valid and the analysis proceeds identically, except that the gradients and subdifferentials should be understood as being defined relative to $\operatorname{aff}(\mathcal{X})$.}}}

    \item The mapping $c = (c_1, \dots, c_m)$ satisfies that each component $c_j : \mathbb{R}^{n_2} \to \mathbb{R}$ is $L_{c_j}$-smooth and twice continuously differentiable on $\cY$. 
     
    \item The function $\bF$ satisfies the following local Kurdyka-\L{}ojasiewicz (KL) condition with respect to $(y, z)$: there exist constants $C > 0$, $\theta \in [1/2, 1)$, $\gamma > 0$, and $\sigma \geq 0$ such that for any {$x \in \cX$},
    \beq \label{poten_KL}
    C (\bF^*(x) - \bF(x, y, z))^{\theta} \leq \dist (0, \partial_{(y, z)} \bF(x, y, z)) \qquad \forall (y, z) \in \overline{\cL}(x),
    \eeq
    where
    \beq \label{poten_Lev}
    \overline{\cL}(x)
    := \{(y, z): 0 < \bF^*(x) - \bF(x, y, z)
         \leq \gamma\,\dist(0,\partial\bPsi(x))^\sigma\}.
    \eeq
\end{enumerate}
\end{assumption}

We now make some remarks on Assumption \ref{ass:Lip_Smo_KL}.

\begin{remark}
\begin{itemize}
\item[(i)] {The Lipschitz continuity of $f(\cdot, y)$ in Assumption~\ref{ass:Lip_Smo_KL}(i) is used to establish the Lipschitz continuity of the maximal function $F^*$. Consequently, by Rademacher's theorem, it follows that $F^*$ is differentiable almost everywhere. This property plays a key role in establishing the local generalized H\"{o}lder smoothness property of $F^*$ in Theorem~\ref{thm:Phi_localHolder}.}


\item[(ii)] {The twice continuous differentiability of the mapping $c$ in Assumption~\ref{ass:Lip_Smo_KL}(iii) is imposed to ensure that the constraint function $\bc$ in the lifted minimax problem~\eqref{lift-minimax} is continuously differentiable. This property will be used to analyze Algorithm~\ref{algo_implement:subsolver} for solving $\max_y F(x,y)$.}


\item[(iii)] {Assumption~\ref{ass:Lip_Smo_KL}(iv) requires that $\bF(x,\cdot,\cdot)$ satisfy the KL property only on the level set $\overline{\cL}(x)$, which we refer to as a local KL condition (in contrast to its global counterpart). When $\sigma > 0$, this level set depends on $x$, in particular on $\dist(0,\partial \bPsi(x))$. If $x$ is far from a stationary point of $\bPsi$, then $\overline{\cL}(x)$ is relatively large, and hence the KL property holds on a correspondingly large level set of the maximization problem in~\eqref{lift-minimax}. As $x$ approaches a stationary point, however, the level set shrinks, and the KL property is required only on a correspondingly smaller level set. Moreover, this local KL condition is used to analyze Algorithm~\ref{algo_implement:subsolver}. Since this algorithm is based on the reformulation~\eqref{equiv-prob} of problem~\eqref{intro:problem}, it is natural to impose the local KL condition on $\bF(x,\cdot,\cdot)$.}

\item[(iv)] For each $x \in \cX$, if $\bF(x,\cdot,\cdot)$ is semi-algebraic, subanalytic, or a structured nonsmooth function, then it satisfies a local KL condition with its own KL constant, KL exponent, and corresponding level set (see, e.g., \cite{attouch2010proximal, bolte2007lojasiewicz, bolte2007clarke}). Consequently, for such $\bF(x,\cdot,\cdot)$, if there exist a common KL constant and exponent for all $x \in \cX$, and if the level set does not shrink faster than $\mathcal{O}(\dist(0,\partial \bPsi(x))^\sigma)$ for some $\sigma > 0$ as $x$ approaches a stationary point, then Assumption~\ref{ass:Lip_Smo_KL}(iv) is satisfied. Otherwise, it may be challenging to develop a first-order method with attractive complexity guarantees for finding an approximate stationary point of problem~\eqref{intro:problem}.

\end{itemize}
\end{remark}

We end this subsection with a simple constrained minimax problem that satisfies the local KL condition in Assumption~\ref{ass:Lip_Smo_KL}(iv) but not the global one (which requires the KL inequality to hold for all $(y, z) \in \dom\,\bF(x, \cdot, \cdot)$), thereby illustrating that the local KL condition applies to a broader class of minimax problems.

\medskip

{\bf Example 1.}
Consider the problem
\[
\min_{1 \leq x \leq 2}\, \max_{c(y) \leq 0} \{ -x^2 \sin y - \delta_{[0, 2\pi/3]}(y) \},
\]
where $c(y) = y^3/6 + y - \pi$, which is $L_c$-smooth on $[0, 2\pi/3]$ with $L_c = 2\pi/3$. 
Notice that this problem is a special case of \eqref{intro:problem} with $f(x, y) = -x^2 \sin y$, $p(x) = \delta_{[1,2]}(x)$, and $q(y) = \delta_{[0, 2\pi/3]}(y)$. We can show that the function $\bF$ defined in \eqref{intro_potential} does not satisfy the global KL condition at any $x \in \dom\,p$, since the KL inequality fails when $y = z = \pi/2$. Indeed, one can observe from the definitions of $\bc$ and $c$ that 
$\bc(\pi/2, \pi/2) = c(\pi/2) = \pi^3/48 - \pi/2 < 0$, 
which along with \eqref{intro_potential} yields 
$\partial_{(y, z)} \bF(x, \pi/2, \pi/2) = (-x^2\cos(\pi/2), 0) = (0,0)$.
On the other hand, note that $\bF^*(x)=0$ for all $x \in \dom\,p$. Using this,  \eqref{intro_potential}, \eqref{intro_PhiOver}, and $\bc(\pi/2, \pi/2) \leq 0$, one has
\[
\bF^*(x) - \bF(x, \pi/2, \pi/2) = 0 - (-x^2 \sin(\pi/2)) = x^2 \geq 1 \qquad \forall x \in \dom\,p.
\]
Combining these, we see that the KL inequality does not hold for $y = z = \pi/2$ at any $x \in \dom\,p$.

However, it can be shown that the KL condition holds on the level set $\overline{\cL}(x)$ defined in \eqref{poten_Lev} with $\gamma=1/2$ and $\sigma=0$ for all $x \in \dom\,p$. To this end, it suffices to show that for any $x \in \dom\,p$, the condition holds on a larger level set $\widehat{\cL}_x = \{(y, z): 0 < \bF^*(x) - \bF(x, y, z) \leq x^2/2\}$. Indeed, one can show that $\widehat{\cL}_x = \{(y, z): y \in (0, \pi/6],\, z \in [0, 2\pi/3],\, \bc(y, z) \leq 0\}$ for all $x \in \dom\,p$. Moreover, it can be verified that $\bc(y, z) < 0$ for all $(y, z) \in (0, \pi/6] \times [0, 2\pi/3]$, due to $\bc(0, z) < 0$, $\bc(\pi/6, z) < 0$ for all $z \in [0, 2\pi/3]$, and  the convexity of $\bc(\cdot, z)$. In view of these and \eqref{intro_potential}, we see that $\partial_{(y, z)} \bF(x, y, z) = (-x^2 \cos y, -\cN_{[0, 2\pi/3]}(z))$ for all $x \in \dom\,p$ and $(y, z) \in \widehat{\cL}_x$. Then, one can verify that
\[
C (\bF^*(x) - \bF(x, y, z))^{\theta} \leq \dist (0, \partial_{(y, z)} \bF(x, y, z)) \qquad \forall (y, z) \in \widehat{\cL}_x
\]
for all $x \in \dom\,p$ with $C = (3/2)^{1/2}$ and $\theta = 1/2$.

\section{Theoretical properties of problem \eqref{intro:problem}} \label{sec:localKL_Holder}

In this section, we establish several theoretical properties for problem \eqref{intro:problem}, which will be used later for algorithm design and analysis.

The following result shows that $F^*$ possesses a local {generalized} H\"{o}lder smoothness property. Its proof is deferred to Subsection~\ref{sec:proof1}.

\begin{theorem} \label{thm:Phi_localHolder}
Suppose that Assumption \ref{ass:Lip_Smo_KL} holds.  Let $\epsilon>0$ be given and 
\beq \label{U-eps}
\cU_\epsilon := \{x \in \cX: \dist(0, \partial \Psi(x)) > \epsilon\}.
\eeq 
Then the following statements hold.
\begin{itemize}
\item[(i)] $\partial^\rC_\cX F^*(x)$ is a singleton for all $x \in \cU_\epsilon$, and $F^*$ is differentiable on $\cU_\epsilon \cap \operatorname{int}(\cX)$. 
\item[(ii)]  For any $x, x' \in \cU_\epsilon \cap \operatorname{int}(\cX)$ satisfying $\|x - x'\| \leq \gamma \epsilon^\sigma/(2L_f)$, we have
\[
\|\nabla F^*(x) - \nabla F^*(x')\| \leq L_{\nabla f}\|x - x'\| + (1-\theta)^{-1} {C^{-1/\theta}}L_{\nabla f}^{1/\theta}\, \|x- x'\|^{\frac{1-\theta}{\theta}}.    
\]
\item[(iii)] For any $x, x' \in \cU_\epsilon$ satisfying $\|x - x'\| \leq \gamma \epsilon^\sigma/(4L_f)$, we have
\[
\|\nabla^\rC_\cX F^*(x) - \nabla^\rC_\cX F^*(x')\| \leq L_{\nabla f}\|x - x'\| + (1-\theta)^{-1} {C^{-1/\theta}}L_{\nabla f}^{1/\theta}\, \|x- x'\|^{\frac{1-\theta}{\theta}}.
\]
\item[(iv)] It holds that
\[\nabla^\rC_\cX F^*(x) = \nabla_x f(x, y^*) \qquad \forall x\in\cU_\epsilon, \ y^* \in Y^*(x).\]
\end{itemize}
\end{theorem}    

{The following result is a consequence of Theorem~\ref{thm:Phi_localHolder}, whose proof is identical to that of \cite[Corollary 1]{lu2025first} and thus omitted.}

\begin{corollary} \label{cor:Fstar-est}
Let $\epsilon>0$ be given and $\cU_\epsilon$ be defined in \eqref{U-eps}. Suppose that Assumption~\ref{ass:Lip_Smo_KL} holds. Then, for any $x, x'$ satisfying $[x, x'] \subseteq \cU_\epsilon$ and $\|x - x'\| \leq \gamma \epsilon^\sigma/(4L_f)$, we have
\beq \label{cor1_Fstar_est}
F^*(x) \leq F^*(x') + \langle \nabla^\rC_\cX F^*(x'), x - x' \rangle + \frac{1}{2} L_{\nabla f} \|x - x'\|^2 + \frac{M}{1+\nu} \|x - x'\|^{1+\nu},
\eeq
where 
\beq \label{Phi_smooth_constants}
M := (1-\theta)^{-1} C^{-1/\theta} L_{\nabla f}^{1/\theta}, \quad \nu := \theta^{-1}(1-\theta).
\eeq
\end{corollary}    

The theorem below establishes a local $(1-\theta)^{-1}$-growth property of $F(x, \cdot)$ for every $x \in \cX$, whose proof is deferred to Subsection~\ref{sec:proof1}.

\begin{theorem} \label{thm:Phi_QG}
Suppose that Assumption~\ref{ass:Lip_Smo_KL} holds. Then it holds that for any $x \in \cX$,  
\beq \label{lem5_Phi_QG}
F^*(x) - F(x, y) \geq (C(1-\theta))^{\frac{1}{1-\theta}} \dist(y, Y^*(x))^{\frac{1}{1-\theta}} \qquad \forall y \in \cL(x),
\eeq
where
\beq \label{Lev_set}
\cL(x) := \{y : 0< F^*(x)-F(x, y) \leq \gamma\, \dist(0, \partial \Psi(x))^\sigma\}.
\eeq
\end{theorem}

We next introduce the Mangasarian–Fromovitz constraint qualification (MFCQ) for problem~\eqref{intro:problem}, which will be used  frequently in the paper.

\begin{assumption}[{\bf MFCQ}] \label{ass:MFCQ}
For every $y \in \cS$, there exists some $\tilde{y} \in \cY$ such that
\[
\langle \nabla c_j(y), \tilde{y}-y \rangle < 0 \qquad \forall j \in I(y):=\{i: c_i(y) = 0\},
\]
where $\cY$ and $\cS$ are defined in \eqref{intro_domain}.
\end{assumption}

Under the above MFCQ condition, we establish a key property of the mapping $\bc(\cdot,\cdot)$, which will be used subsequently in the paper.

\begin{theorem} \label{thm:UniSlater} 
Let $\bc$ be defined in \eqref{intro_cOver}. Suppose that Assumptions \ref{ass:Lip_Smo_KL} and \ref{ass:MFCQ} hold. Then there exists a constant $\zeta > 0$ such that for every $z \in \cS$, one can find a point $y \in \cY$ (possibly depending on $z$) satisfying $\bc(y, z) \leq -\zeta$.
\end{theorem}

\begin{proof}
Suppose for contradiction that the conclusion of this theorem does not hold. Then there exist a positive sequence $\{\zeta_k\}$ and a sequence $\{z^k\} \subset \cS$ such that $\lim_{k\to\infty }\zeta_k = 0$ and  
\beq \label{bc-ineq1}
\bc(y, z^k) \nleq -\zeta_k \quad  \forall y\in\cY, \, \forall k.
\eeq
Since $\cY$ is compact and $c$ is continuous, it follows that $\cS$ is compact. Hence, the sequence $\{z^k\}$ admits a convergent subsequence. Passing to such a subsequence if necessary, we may assume without loss of generality that $z^k \to z$ for some $z \in \cS$. By this fact and Assumption \ref{ass:MFCQ}, there exists some $\tz \in \cY$ such that
\[
\langle \nabla c_{j}(z), \tz-z \rangle < 0 \qquad \forall j \in I(z) = \{i: c_i(z) = 0\}.
\]
Let $y(t)=z+t(\tz-z)$ for all $t$. Then, for each $j \in I(z)$, one has
\[
\lim_{t \downarrow 0} \frac{\bc_{j}(y(t), z)}{t}= \langle \nabla c_{j}(z), \tz-z \rangle < 0,
\]
which implies that $\bc_{j}(y(t), z) < 0$ for all sufficiently small $t > 0$. In addition, for each $j \notin I(z)$,  we have $\bc_{j}(y(t), z) < 0$ for all sufficiently small $t > 0$, thanks to $\bc_{j}(z, z) = c_j(z) < 0$ and the continuity of $\bc_j$. Moreover, it follows from $z, \tz \in\cY$ and the convexity of $\cY$ that $y(t)\in\cY$ for all $t\in [0,1]$. Consequently, there exists some $\bar{y} \in \cY$ such that 
\beq \label{bc-ineq2}
\bc(\bar{y}, z) < 0.
\eeq

On the other hand, since $\bar{y} \in \cY$, it follows from \eqref{bc-ineq1} that $\bc(\bar{y}, z^k) \nleq -\zeta_k$ for all $k$. Hence, there exists some $\bar{j}$ such that $\bc_{\bar{j}}(\bar{y}, z^k) > -\zeta_k$ holds for infinitely many $k$. Passing to a subsequence if necessary, we may assume without loss of generality that this inequality holds for all $k$. Taking limits on both sides of this inequality as $k \to \infty$ and using $\lim_{k\to\infty }\zeta_k = 0$, $\lim_{k\to\infty}z^k = z$, and the continuity of $\bc_{\bar{j}}$, we obtain that $\bc_{\bar{j}}(\bar{y}, z) \geq 0$, which contradicts \eqref{bc-ineq2}. Hence, the conclusion of this theorem holds.
\end{proof}

\section{A sequential convex programming method for constrained KL function minimization} \label{sec:subsolver_SCP}

In this section, we consider constrained optimization problems of the form
\beq \label{pr:KL_opt}
    h^* = \min_{c(z) \leq 0} \{ h(z) := g(z) + q(z) \},
\eeq
where $q$ and $c$ are defined in Section \ref{sec:intro} that satisfy Assumption \ref{ass:Lip_Smo_KL}, and $g$ and $h$ satisfy the following assumption.

\begin{assumption} \label{assum:g-h}
The function $g$ is $L$-smooth on $\cY$, and $\bh$ satisfies the following KL condition:
\beq \label{KL_subpr}
 C(\bh(z, w) - \bh^*)^{\theta} \leq \dist(0, \partial \bh(z, w)) \qquad \forall (z, w) \text{ with } \bh^* < \bh(z, w) \leq \bh^* + \eta
\eeq
for some constants $C,\eta>0$ and $\theta \in [1/2, 1)$,  where $\cY=\dom\,q$, and  
 \beq \label{subpr_overh}
\bh(z, w) := h(z) + \delta_{\bc(\cdot,\cdot) \leq 0}(z, w) + \delta_{\cY}(w), \qquad \bh^* := \min_{z, w} \bh(z, w),
\eeq
and $\bc$ is given in \eqref{intro_cOver}. 
\end{assumption}

Sequential convex programming (SCP) methods have been studied in the literature (see, e.g., \cite{lu2012sequential, yu2021convergence}) for solving problem~\eqref{pr:KL_opt}. These methods solve a sequence of simple convex optimization problems of the form~\eqref{surro_argmin}. In particular, \cite{lu2012sequential, yu2021convergence} apply line search schemes to both the objective and the constraints of~\eqref{pr:KL_opt} to generate a sequence of feasible points that ensures sufficient reduction in the objective function $h$ along the sequence. Motivated by these works, we propose a variant of the SCP method for solving~\eqref{pr:KL_opt}, which will subsequently serve as a subroutine for solving problem~\eqref{intro:problem}. To suit our purpose, we apply the line search scheme only to the objective, so that the resulting sequence is stronger than those generated by the SCP methods in \cite{lu2012sequential, yu2021convergence}---in particular, it satisfies the constraints more strictly. Specifically, at each iteration, the method performs multiple constrained proximal gradient steps using the surrogate constraint mapping $\bc$, along with a backtracking line search to ensure sufficient reduction in $h$. The method terminates once a practical stopping criterion---designed to guarantee that $\dist(0, \partial \bh(z^{k+1}, z^k))$ is sufficiently small---is met for some $z^k$ and $z^{k+1}$. The proposed method is described in Algorithm~\ref{algo_implement:subsolver}, where $\bc$ is defined in~\eqref{intro_cOver}.

\begin{algorithm}[H]
\caption{A sequential convex programming method for problem~\eqref{pr:KL_opt}}\label{algo_implement:subsolver}
\begin{algorithmic}[1]
\REQUIRE $\{L_{c_j}\}_{j=1}^m$ from Assumption~\ref{ass:Lip_Smo_KL}; $\underline{L} > 0$, $\rho > 1$, $\beta > 0$, $\tau > 0$, and a point $z^0 \in \{ z : h(z) \leq h^* + \eta,\; c(z) \leq 0 \}$.
\FOR{$k=0,1,2,\dots$}
    \FOR{$i=0,1,2,\dots$}
        \STATE Set $L_{k,i} = \underline{L} \rho^i$.
        \STATE Compute
        \beq \label{surro_argmin}
        \begin{aligned}
            z^{k+1, i} = &\mathop{\arg\min}_z \Big\{ \langle \nabla g(z^k), z \rangle + \frac{L_{k, i}}{2}\|z - z^k\|^2 + q(z) \Big\} \\
            &\qquad \text{s.t.} \quad \bc(z, z^k) \leq 0,
        \end{aligned}
        \eeq
        and its optimal Lagrange multiplier $\lambda^{k, i}$.
        \IF{$h(z^{k+1,i}) \leq h(z^k)- \beta\|z^{k+1,i} - z^k\|^2/2$}
            \STATE Set $z^{k+1} = z^{k+1,i}$, $L_k = L_{k, i}$, $\lambda^k = \lambda^{k, i}$.
            \STATE \textbf{break}
        \ENDIF
    \ENDFOR
    \STATE Terminate the algorithm and output $z^{k+1}$ if
    \beq \label{subsolver_terminate}
    \|\nabla g(z^{k+1}) - \nabla g(z^k) - L_k(z^{k+1}-z^k)\|^2 + 4 \Big(\sum_{j=1}^m \lambda^k_j L_{c_j} \Big)^2 \|z^{k+1}-z^k\|^2 \leq \tau^2.
    \eeq
\ENDFOR
\end{algorithmic}
\end{algorithm}

We now make some remarks on subproblem~\eqref{surro_argmin} in Algorithm~\ref{algo_implement:subsolver}. By rearranging the terms in the constraint functions of~\eqref{surro_argmin}, one can see that the constraint is equivalent to
\[
z \in \bigcap_{j=1}^{m} \cB\big( s^{k, j}, \sqrt{R_{k, j}} \big),
\]
where
\[
s^{k, j} = z^k - \nabla c_j(z^k)/L_{c_j}, \quad R_{k, j} = \| \nabla c_j(z^k) \|^2 / L_{c_j}^2 - 2c_j(z^k)/L_{c_j}.
\]
Thus, the constraint in~\eqref{surro_argmin} corresponds to the intersection of Euclidean balls. 
{In addition, since the function $q$ is a simple closed convex function, subproblem~\eqref{surro_argmin}  can typically be reformulated as a convex conic optimization problem, and a customized primal-dual interior-point method (IPM) can be applied to compute its optimal solution and the associated Lagrange multipliers. Moreover, as the Hessians of the objective and the constraint functions in~\eqref{surro_argmin} are multiples of the identity matrix, the Newton system arising in the IPM can be solved cheaply when the epigraph of $q$ is polyhedral (e.g., $q(\cdot)=\|\cdot\|_1$).}

The following theorem establishes the well-definedness of Algorithm~\ref{algo_implement:subsolver} and several key relations used in the subsequent analysis. Similar results were obtained in \cite[Lemma 2.4]{yu2021convergence} for a more sophisticated SCP method. Since our SCP method is simpler, we provide a more concise proof in Subsection~\ref{sec:proof2} for reference.

\begin{theorem} \label{thm:sub_well} 
Suppose that Assumptions \ref{ass:Lip_Smo_KL}, \ref{ass:MFCQ}, and \ref{assum:g-h} hold. 
Let $\{L_{c_j}\}_{j=1}^m, L, \underline{L}, \rho, \beta$ be given in Assumptions \ref{ass:Lip_Smo_KL} and \ref{assum:g-h}, and Algorithm~\ref{algo_implement:subsolver}, respectively, $\bar{i} = \lceil \log_\rho ((\beta+L)/(2\underline{L})) \rceil_+$, and $\{L_k\}$, $\{z^k\}$,  and $\{\lambda^k\}$ be generated in Algorithm~\ref{algo_implement:subsolver}. Then the following statements hold.
\begin{enumerate} [label=(\roman*)]
    \item The subproblem \eqref{surro_argmin} has an optimal solution $z^{k+1, i}$ and an optimal Lagrange multiplier $\lambda^{k, i}$, and the inner loop terminates in at most $\bar{i}+1$ iterations and outputs a point $z^{k+1} \in \dom\,q$ satisfying $c(z^{k+1}) \leq 0$ at each outer iteration $k$.
    \item For each $k$, it holds that
    \begin{align} 
   & \underline{L} \leq L_k \leq \bar{L} := \max\{ \underline{L}, (\beta+L)\rho/2\}, \label{Lk-bound} \\
    & \lambda^k \geq 0, \quad \lambda^{k}_j \Big(c_j(z^k) + \langle \nabla c_j(z^k), z^{k+1}-z^k \rangle + \frac{L_{c_j}}{2}\|z^{k+1}-z^k\|^2\Big) = 0 \quad \forall j \in \{1,...,m\}, \label{thm2_KKT_CompSlack} \\
    & 0 \in \nabla g(z^k) + \Big(L_{k} + \sum_{j=1}^m \lambda^{k}_j L_{c_j}\Big)(z^{k+1}-z^k) + \partial q(z^{k+1}) + \sum_{j=1}^m \lambda^{k}_j \nabla c_j(z^k). \label{thm2_KKT_Stationary}
   \end{align}
\end{enumerate}
\end{theorem}

To establish the convergence rate of Algorithm~\ref{algo_implement:subsolver}, we now present a result that provides an upper bound on the sequence of Lagrange multipliers $\{\lambda^k\}$ generated by Algorithm~\ref{algo_implement:subsolver}. Its proof is deferred to Subsection~\ref{sec:proof2}.

\begin{lemma} \label{lambda-bound}
Suppose that Assumptions \ref{ass:Lip_Smo_KL}, \ref{ass:MFCQ}, and \ref{assum:g-h} hold.  Let $\lambda^k$ be generated in the $k$th outer iteration of Algorithm~\ref{algo_implement:subsolver} for some $k\geq 0$, $\zeta$ be given in Theorem~\ref{thm:UniSlater}, $D_{\cY}$, $M_q$, $\bar{L}$ be given in \eqref{DY-Mq} and \eqref{Lk-bound}, and let $G= \max_{z \in \cY } \|\nabla g(z)\|$. Then it holds that
\beq \label{prop1_explicitBound}
\|\lambda^k\|_1 \leq A:=\zeta^{-1}(G D_{\cY} + \bar{L}D_{\cY}^2/2 + M_q).
\eeq
\end{lemma}

The theorem below shows that Algorithm~\ref{algo_implement:subsolver} terminates in a finite number of iterations and outputs a desired approximate solution to problem \eqref{pr:KL_opt}. Its proof is deferred to Subsection~\ref{sec:proof2}.

\begin{theorem} \label{thm:sub_iters}
Suppose that Assumptions \ref{ass:Lip_Smo_KL}, \ref{ass:MFCQ}, and \ref{assum:g-h} hold. Let $\{L_{c_j}\}_{j=1}^m, L, C, \theta, \eta, \bar{L}, A, \beta, \tau$ be given in Assumptions \ref{ass:Lip_Smo_KL} and \ref{assum:g-h}, \eqref{Lk-bound}, \eqref{prop1_explicitBound}, and Algorithm~\ref{algo_implement:subsolver}, respectively, and let
\begin{align}
&\omega = \Big( (L+\bar{L})^2 + 4A^2 \big(\sum_{j=1}^m L_{c_j}\big)^2 \Big)^{\frac{1}{2}}, \label{lem9_omega} \\
&\alpha = \frac{\beta C^2}{2 \omega^2}, \qquad C' = \min\left\{ \frac{\alpha}{2}, \frac{(2^{\frac{2\theta-1}{2\theta}} - 1)\eta^{1-2\theta}}{2\theta-1}\right\}, \label{lem4:alpha_Cp} \\
&\overline{K}_\theta :=
\begin{cases}
\Big\lceil \log_{1+\alpha}\big(\frac{2\omega^2\eta}{\beta \tau^2}\big) \Big\rceil_{+} + 1  \quad & \text{if } \theta = \tfrac{1}{2}, \\[2ex]
\Big\lceil \frac{1}{C' (2\theta-1)} \big( \frac{2\omega^2}{\beta\tau^2} \big)^{2\theta-1} \Big\rceil + 1 \quad & \text{if } \theta \in (\tfrac{1}{2}, 1). \label{thm3:overK}
\end{cases}  
\end{align}
Then Algorithm~\ref{algo_implement:subsolver} terminates in at most $\overline{K}_\theta$ outer iterations and outputs a point $z^{k+1}$ satisfying 
\beq \label{thm3:output}
h(z^{k+1}) - h^* \leq (C^{-1}\tau)^{\frac{1}{\theta}}
\eeq
for some $k < \overline{K}_\theta$.
\end{theorem}

\section{An inexact proximal gradient method for problem \eqref{intro:problem}} \label{sec:FOD}

In this section, we propose an inexact proximal gradient method for solving problem \eqref{intro:problem} and analyze its complexity for finding an $(\gamma \epsilon^\sigma/(4L_f),\epsilon)$-stationary point of \eqref{intro:problem} for $\epsilon >0$. 

Before proceeding, we introduce some additional notation below. Given any $\epsilon>0$, let
\begin{align}
&\cX_\epsilon := \{x\!\in\!\cX : \dist(0, \partial\Psi(x)) \leq \epsilon\}, \quad \cX^\rc_\epsilon :=  \{x\!\in\!\cX : \dist(x, \cX_\epsilon ) > \gamma \epsilon^\sigma / (4L_f) \}, \quad r:=\gamma \epsilon^\sigma / (4L_f)  \label{X-eps-r}
\end{align}
where $\gamma, \sigma, L_f$ are given in Assumption \ref{ass:Lip_Smo_KL}.

To propose a method for finding an $(r,\epsilon)$-stationary point of problem \eqref{intro:problem}, we first make several key observations. Suppose $x'\in\cX^\rc_\epsilon$, that is, $x'$ is not an $(r,\epsilon)$-stationary point of \eqref{intro:problem}. For any $x\in\cX\cap\cB(x',r)$, we observe that $[x',x] \subseteq \cX$ and moreover $\dist(0,\partial\Psi(z))>\epsilon$ for all $z\in [x',x]$. In view of these observations and the definition of $\cU_\epsilon$ in \eqref{U-eps}, it follows that $[x', x] \subseteq \cU_\epsilon$. {Hence, by Corollary \ref{cor:Fstar-est},  the relation \eqref{cor1_Fstar_est} holds for such $x$ and $x'$. In addition, note from $\theta \in [1/2, 1)$ and \eqref{Phi_smooth_constants} that $\nu \in (0, 1]$. By this relation and  \cite[(2.15)]{nesterov2015universal} with $M_\nu$ and $t$ replaced by $M$ and $ \|x - x'\|$, respectively, one has
 \[
 M(1+\nu)^{-1} \|x - x'\|^{1+\nu} \leq \frac12 \left[\frac{1-\nu}{1+\nu} \cdot \frac{1}{\delta}\right]^{\frac{1-\nu}{1+\nu}}M^{\frac{2}{1 + \nu}} \|x - x'\|^2+\frac{\delta}{2} \leq \frac12\Big( \delta^{\frac{\nu - 1}{1 + \nu}} M^{\frac{2}{1 + \nu}} \|x - x'\|^2 + \delta\Big) \qquad \forall \delta>0,
 \]
where the first relation follows from Young's inequality (see \cite[(2.15)]{nesterov2015universal}), and the second relation is due to $\nu \in (0, 1]$. Combining this inequality with \eqref{cor1_Fstar_est}, we obtain that
\beq  \label{Fstar-bound} 
F^*(x) \leq F^*(x') + \langle {\nabla^\rC_\cX F^*(x')}, x \!-\! x' \rangle + \frac{1}{2} \big(L_{\nabla f} \!+\! \delta^{\frac{\nu - 1}{1 + \nu}} M^{\frac{2}{1 + \nu}}\big) \|x \!-\! x'\|^2 + \frac{\delta}{2} \qquad \forall x\in\cX\cap \cB(x',r).
\eeq
By this relation and $\Psi(\cdot)=F^*(\cdot)+p(\cdot)$, we further have
\[
\Psi(x) \leq F^*(x') + \langle {\nabla^\rC_\cX F^*(x')}, x - x' \rangle + \frac{1}{2} \big(L_{\nabla f} + \delta^{\frac{\nu - 1}{1 + \nu}} M^{\frac{2}{1 + \nu}}\big) \|x - x'\|^2 +p(x)+ \frac{\delta}{2} \qquad \forall x\in\cX\cap \cB(x',r). 
\]}
Therefore, when $x'\in\cX$ is not an $(r,\epsilon)$-stationary point of \eqref{intro:problem}, $\Psi$ is bounded above by a much simpler function that is the sum of a simple quadratic function and $p(\cdot)$ in a neighborhood of $x'$.

Based on the above observation, it is natural to propose a proximal-gradient-type method for finding an $(r,\epsilon)$-stationary point of problem \eqref{intro:problem}. Specifically, the method generates the sequence $\{x^k\}$ according to
\[
x^{k+1}=\mathop{\arg\min}_{x \in \cB(x^k, r)} \Big\{\langle {\nabla^\rC_\cX F^*(x^k)}, x  \rangle + \frac{1}{2} \bar{L}_k \|x - x^k\|^2 +p(x)\Big\}
\]
with $\bar{L}_k= L_{\nabla f} + \delta_k^{(\nu-1)/(1+\nu)}M^{2/(1+\nu)}$ for a suitable choice of $\delta_k>0$. However, since $F^*$ is a maximal function, the exact value of ${\nabla^\rC_\cX F^*(x^k)}$ is generally unavailable. To overcome this difficulty, we approximate ${\nabla^\rC_\cX F^*(x^k)}$ by $\nabla_x f(x^k, y^k)$, where $y^k$ is a suitably chosen approximate solution of the $k$th subproblem
\beq \label{subprob-k}
\min_y\{-f(x^k, y) + q(y): c(y) \leq 0\}. 
\eeq
Such $y^k$ is obtained using Algorithm~\ref{algo_implement:subsolver}, initialized from $y^{k-1}$ (see line 5 of Algorithm~\ref{algo_implement:whole}). We show that if $y^0$ is a suitable approximate solution to the initial subproblem and $\{x^\ell\}_{0 \leq \ell < k}$ are not $(\gamma \epsilon^\sigma/(4L_f),\epsilon)$-stationary points of \eqref{intro:problem}, then $y^k$ generated in this manner is indeed a desired approximate solution to \eqref{subprob-k} (see Lemma~\ref{lem:subsolver_output}). 

We are now ready to present an inexact proximal gradient method for solving problem \eqref{intro:problem}.

\begin{algorithm}[H]
\caption{An inexact proximal gradient method for problem \eqref{intro:problem}}\label{algo_implement:whole}
\begin{algorithmic}[1]
\REQUIRE $L_f$, $L_{\nabla f}$, $\{L_{c_j}\}_{j=1}^m$, $C$, $\theta$, $\gamma$, $\sigma$ from Assumption~\ref{ass:Lip_Smo_KL}; $\epsilon > 0$, $\underline{L} > 0$, $\rho > 1$, $\beta > 0$, and initial point $(x^0, y^0) \in \cX \times \cY$ satisfying $F^*(x^0) - F(x^0, y^0) \leq \min\{\gamma\epsilon^\sigma/2, 1\}$.

\STATE Set $r = \gamma \epsilon^\sigma / (4L_f)$, $M = (1-\theta)^{-1} C^{-1/\theta} L_{\nabla f}^{1/\theta}$, $\nu = \theta^{-1}(1-\theta)$.
\FOR{$k=0,1,2,\dots$}
    \STATE Set $\delta_k = 1/(k+1)$, $\eta_{k} = 1/(k+1)$,  $\bar{L}_k= L_{\nabla f} + \delta_k^{(\nu-1)/(1+\nu)}M^{2/(1+\nu)}$.
    \STATE Compute 
    \[
        x^{k+1} = \mathop{\arg\min}_{x \in \cB(x^k, r)} \Big\{ \langle \nabla_x f(x^{k}, y^{k}), x \rangle + \frac{\bar{L}_k}{2} \|x - x^k\|^2 + p(x) \Big\}.
    \]
    \STATE Call Algorithm~\ref{algo_implement:subsolver} with $g(\cdot) \gets -f(x^{k+1}, \cdot)$, $q(\cdot) \gets q(\cdot)$, $c(\cdot) \gets c(\cdot)$, $z^0 \gets y^k$, $\underline{L} \gets \underline{L}$, $\rho \gets \rho$, $\beta \gets \beta$, $\{L_{c_j}\}_{j=1}^m \gets \{L_{c_j}\}_{j=1}^m$, $\tau \gets C \min\Big\{ (\tfrac{1}{2} \gamma \epsilon^\sigma)^{\theta}, \eta_{k+1}^{\frac{\theta}{2(1-\theta)}} \Big\}$, and denote its output as $y^{k+1}$.
\ENDFOR
\end{algorithmic}
\end{algorithm}

\begin{remark}
(i) {The required $y^0$ for Algorithm~2 can be obtained by applying Algorithm~1 to the problem $\max_y F(x^0,y)$, provided that an initial point in a region where the local KL condition holds is readily available. Alternatively, such a $y^0$ can be efficiently computed when $F(x^0,\cdot)$ is concave. Moreover, even when $F(x^0,\cdot)$ is nonconcave, the problem corresponding to the specific choice of $x^0$ may still possess a favorable structure, making it computationally tractable to compute such a $y^0$ by exploiting this structure.}

(ii) Some input parameters required by Algorithm~\ref{algo_implement:whole} may not be readily available in practice. It would therefore be worthwhile to develop a parameter-free variant of Algorithm~\ref{algo_implement:whole}. Alternatively, in practical implementations, one may run the algorithm with a range of trial parameters and adjust them until its performance stabilizes.
\end{remark}

To analyze the complexity of Algorithm \ref{algo_implement:whole} for computing an $(\gamma \epsilon^\sigma/(4L_f),\epsilon)$-stationary point of problem~\eqref{intro:problem}, it is necessary to establish that the sequence of Lagrange multipliers generated by Algorithm~\ref{algo_implement:subsolver} for solving subproblem~\eqref{subprob-k} remains uniformly bounded, independent of $k$.

\begin{lemma} \label{lambda-uniform-bound}
Suppose that Assumption \ref{ass:Lip_Smo_KL} holds.  Let $\{\lambda^{k,\ell}\}$ denote the sequence of Lagrange multipliers generated by Algorithm~\ref{algo_implement:subsolver} during the $k$th iteration of Algorithm~\ref{algo_implement:whole}, and let $\zeta$ be given in Theorem~\ref{thm:UniSlater}, $\bar{L}$ be given in \eqref{Lk-bound}, $D_\cY, M_q$ be defined in \eqref{DY-Mq},  and 
\[
G_f = \max\{\|\nabla_y f(x, y)\|:(x, y) \in\ \cX \times \cY\}.
\]
Suppose that $G_f < \infty$. Then it holds that
\beq \label{Af}
\|\lambda^{k,\ell}\|_1 \leq A_f:=\zeta^{-1}(G_f D_{\cY} + \bar{L}D_{\cY}^2/2 + M_q) \qquad \forall k, \ell.
\eeq
\end{lemma}

The following theorem establishes an \emph{iteration complexity} bound for Algorithm~\ref{algo_implement:whole} to compute an $(\gamma \epsilon^\sigma/(4L_f),\epsilon)$-stationary point of problem~\eqref{intro:problem} for any $\epsilon \in (0, 1/e]$. The proof is deferred to Subsection~\ref{sec:proof3}.

\begin{theorem} \label{thm:outer_bound}
Let $L_f, L_{\nabla f}, C, \theta, \gamma, \sigma$ be given in Assumption~\ref{ass:Lip_Smo_KL}, $M, \nu$ be defined in \eqref{Phi_smooth_constants}, $\epsilon$ be given in Algorithm~\ref{algo_implement:whole}, and
\begin{align*}
& \widehat{A} = (1-\theta)^{-2} C^{-2}L_{\nabla f}^2, \quad \widehat{L} = L_{\nabla f} + M^{2/(1+\nu)}, \\
& a = 8 \big( \Psi(x^0) - \Psi^* + 3 + 2\widehat{L}^{-1}\widehat{A}\,\big), \quad b = 8 \big(3/2 + \widehat{L}^{-1}\widehat{A}\,\big), \\
& \widehat{C}_1 = \Big( 36(1+\nu)\nu^{-1}b\widehat{L} \lceil \log(18 (1+\nu)\nu^{-1}b\widehat{L}) \rceil_+ + 72(1+\nu)\nu^{-1}b\widehat{L} + 1 \Big)^{\frac{1+\nu}{2\nu}}, \\
& \widehat{C}_2 = \Big( \frac{4b(1+\nu)(3M)^{2/\nu}}{M^{2/(1+\nu)}} \Big\lceil\! \log\Big(\frac{2b(1+\nu)(3M)^{2/\nu}}{M^{2/(1+\nu)}}\Big) \!\Big\rceil_+ + \frac{8b(1+\nu)(3M)^{2/\nu}}{\nu M^{2/(1+\nu)}} + 1 \Big)^{\frac{1+\nu}{2}},  \\
& \widehat{C}_3 = \big( 36\widehat{L} a \big)^{\frac{1+\nu}{2\nu}} + M^{-1} \big( 4a(3M)^{2/\nu} \big)^{\frac{1+\nu}{2}}, \quad \widehat{C}_4 = 72\widehat{A}, \\
& \widehat{C}_5 = \Big( \frac{144(1+\nu)bL_{\nabla f}^2}{M^{2/(1+\nu)}} \Big\lceil\! \log\Big(\frac{72(1+\nu)bL_{\nabla f}^2}{M^{2/(1+\nu)}}\Big) \!\Big\rceil_+ + \frac{288(1+\nu)bL_{\nabla f}^2}{M^{2/(1+\nu)}} + 1 \Big)^{\frac{1+\nu}{2}}, \\
& \widehat{C}_6 = (144aL_{\nabla f}^2)^{\frac{1+\nu}{2}} / M, \\
& \widehat{C}_7 = \Big( \frac{64(1+\nu)bL_f^2}{\gamma^2 M^{2/(1+\nu)}} \Big\lceil\! \log\Big(\frac{32(1+\nu)bL_f^2}{\gamma^2 M^{2/(1+\nu)}}\Big) \!\Big\rceil_+ + \frac{128\sigma(1+\nu)bL_f^2}{\gamma^2 M^{2/(1+\nu)}} + 1 \Big)^{\frac{1+\nu}{2}}, \\
& \widehat{C}_8 = (64aL_f^2)^{\frac{1+\nu}{2}} / (\gamma^{1+\nu} M), \\
&\widehat{K}_\epsilon = \Big \lceil \widehat{C}_1 \epsilon^{-\frac{1+\nu}{\nu}} (\log \epsilon^{-1})^{\frac{1+\nu}{2\nu}} + \widehat{C}_2 \epsilon^{-\frac{1+\nu}{\nu}} (\log \epsilon^{-1})^{\frac{1+\nu}{2}} + \widehat{C}_3 \epsilon^{-\frac{1+\nu}{\nu}} + \widehat{C}_4 \epsilon^{-2} \\ 
&\qquad\quad + \widehat{C}_5 \epsilon^{-(1+\nu)} (\log \epsilon^{-1})^{\frac{1+\nu}{2}} + \widehat{C}_6 \epsilon^{-(1+\nu)} + \widehat{C}_7 \epsilon^{-(1+\nu)\sigma} (\log \epsilon^{-1})^{\frac{1+\nu}{2}} + \widehat{C}_8 \epsilon^{-(1+\nu)\sigma} \Big \rceil. 
\end{align*}
Suppose that $\epsilon \in (0, 1/e]$ and Assumptions \ref{ass:Lip_Smo_KL} and \ref{ass:MFCQ} hold. Then Algorithm~\ref{algo_implement:whole} generates a pair $(x^k, y^k)$ in at most $\widehat{K}_\epsilon$ iterations such that $x^k$ is an $(\gamma \epsilon^\sigma/(4L_f),\epsilon)$-stationary point of problem~\eqref{intro:problem} (or equivalently the problem $\min_x \Psi(x)$), and $y^k$ satisfies
\beq \label{yk-opt}
F^*(x^k) - F(x^k, y^k) \leq \min\Big\{\frac{\gamma \epsilon^\sigma}{2}, \frac{1}{k+1} \Big\}, \quad 
\dist\big(y^k, Y^*(x^k)\big) \leq \frac{1}{C(1-\theta)} \min\Big\{\Big(\frac{\gamma}{2}\Big)^{(1-\theta)} \epsilon^{\sigma (1-\theta)}, \frac{1}{\sqrt{k+1}} \Big\}.
\eeq
\end{theorem}

The next result presents a \emph{first-order oracle complexity} bound for Algorithm~\ref{algo_implement:whole}, measured by the number of evaluations of the gradient $\nabla f$, required to generate an $(\gamma \epsilon^\sigma/(4L_f),\epsilon)$-stationary point of problem~\eqref{intro:problem} for any $\epsilon \in (0, 1/e]$. The proof is deferred to Subsection~\ref{sec:proof3}.

\begin{theorem} \label{thm:operations}
Let $\epsilon \in (0, 1/e]$ be given, $\widehat{K}_\epsilon$ be defined in Theorem~\ref{thm:outer_bound}, $L_{\nabla f}, C, \theta, \gamma, \sigma, \{L_{c_j}\}_{j=1}^m$ be given in Assumption~\ref{ass:Lip_Smo_KL}, $M, \nu$ be defined in \eqref{Phi_smooth_constants}, $\underline{L}, \beta, \rho$ be given in Algorithm~\ref{algo_implement:whole}, $A_f$ be given in \eqref{Af}, and let
\begin{align*}
&\bar{L}_{\nabla f} = \max\Big\{ \underline{L}, \frac{(\beta\!+\!L_{\nabla f})\rho}{2} \Big\}, \quad
\omega_f = \Big( (L_{\nabla f}\!+\!\bar{L}_{\nabla f})^2 \!+\! 4 A_f^2 \big(\sum_{j=1}^m L_{c_j}\big)^2  \Big)^{\frac{1}{2}}, \quad \alpha_f = \frac{\beta C^2}{2\omega_f^2}, \\
&C'_f = \min \Big\{ \frac{1}{2}\alpha_f, \frac{(2^{\frac{2\theta-1}{2\theta}}\!-\!1)(\gamma \epsilon^\sigma)^{1-2\theta}}{2\theta-1} \Big\}, \quad
\Lambda = \max\Big\{ ( \tfrac{1}{2} \gamma \epsilon^\sigma )^{-2\theta}, (\widehat{K}_\epsilon+1)^{\frac{\theta}{1-\theta}} \Big\}, \\
&\overline{K}_{f, \theta} =
\begin{cases}
\left\lceil \log_{1+\alpha_f} (2\omega_f^2\beta^{-1}C^{-2}\gamma \epsilon^\sigma \Lambda) \right\rceil_+ + 1 & \text{if } \theta = \tfrac{1}{2}, \\[1ex]
\left\lceil \frac{1}{C_f'(2\theta-1)} \left( 2\omega_f^2\beta^{-1}C^{-2}\Lambda \right)^{2\theta-1} \right\rceil + 1 & \text{if } \theta \in (\tfrac{1}{2}, 1),
\end{cases} \\
&\widehat{N}_\epsilon = \widehat{K}_\epsilon \Big( \Big\lceil \log_\rho \big( \frac{\beta+L_{\nabla f}}{2\underline{L}} \big) \Big\rceil_+ + 1 \Big) \overline{K}_{f, \theta}.
\end{align*}
Suppose that Assumptions \ref{ass:Lip_Smo_KL} and \ref{ass:MFCQ} hold. Then the total number of evaluations of the proximal operators of $p$ and $q$, and the gradient $\nabla f$ performed by Algorithm~\ref{algo_implement:whole} is at most $\widehat{K}_\epsilon$, $\widehat{N}_\epsilon$, and $\widehat{K}_\epsilon + \widehat{N}_\epsilon$, respectively, to generate a pair $(x^k, y^k)$ such that $x^k$ is an $(\gamma \epsilon^\sigma/(4L_f),\epsilon)$-stationary point of problem~\eqref{intro:problem}, and $y^k$ satisfies
\eqref{yk-opt}.
\end{theorem}

\begin{remark} As shown in Theorem \ref{thm:outer_bound}, Algorithm \ref{algo_implement:whole} achieves an iteration complexity of
\[
\mathcal{O}\left(\epsilon^{-\max\left\{\frac{1}{1-\theta},\,\frac{\sigma}{\theta}\right\}} (\log\epsilon^{-1})^{\frac{1}{2(1-\theta)}}\right)
\]
to compute an $(\gamma \epsilon^\sigma/(4L_f),\epsilon)$-stationary point of problem \eqref{intro:problem}. Furthermore, as established in Theorem~\ref{thm:operations}, the algorithm requires $\mathcal{O}\Big( \epsilon^{-\max\left\{\frac{1}{1-\theta},\,\frac{\sigma}{\theta}\right\}} (\log\epsilon^{-1})^{\frac{1}{2(1-\theta)}} \Big)$ evaluations of the proximal operator of $p$, and the following number of evaluations of the proximal operator of $q$ and the gradient $\nabla f$ to compute such an approximate stationary point of~\eqref{intro:problem}:
\[
\left\{\begin{array}{ll}
    \mathcal{O}\left( \epsilon^{-2\max\left\{1,\, \sigma\right\}} (\log \epsilon^{-1})^2 \right)  & \mbox{if } \theta = \tfrac{1}{2}, \\ [8pt]  
     \mathcal{O}\left( \epsilon^{-\frac{2\theta^2 - 2\theta + 1}{1-\theta} \max\left\{\frac{1}{1-\theta},\,\frac{\sigma}{\theta}\right\}} (\log \epsilon^{-1})^{\frac{2\theta^2 - 2\theta + 1}{2(1-\theta)^2}} \right) & \mbox{if } \theta \in (\tfrac{1}{2},1).
\end{array}
\right.
\]
\end{remark}

\section{Numerical results} \label{sec:Numerical}
In this section, we conduct preliminary experiments to evaluate the performance of our proposed method (Algorithm~\ref{algo_implement:whole}).

Consider the following constrained minimax optimization problem:
\begin{equation} \label{numerical_example}
\min_{x} \max_{c(y) \leq 0} \big\{-\|(y+Ax)\odot(y+Bx)\|^2  + 0.01\,\|x - u\|^2 + 0.01\|x\|_1 + \delta_{\cB(0, 2)}(x) - 0.1\|y\|_1 - \delta_{[-2, 2]^{n_2}}(y) \big\},
\end{equation}
where $A, B \in \mathbb{R}^{n_2 \times n_1}$, $u \in \mathbb{R}^{n_1}$, and $\odot$ denotes the Hadamard (elementwise) product. The mapping $c$ is defined as follows.  Assuming $n_2$ is a multiple of $10$, we set the number of constraints as $m = n_2/10$, and for each $j \in \{1,...,m\}$, the $j$th component of $c$ is 
\beq \label{numerical_constraints}
c_j(y) = e^{y_{10j - 9}} + e^{y_{10j - 8}} + \cdots + e^{y_{10j}} - 10.
\eeq

For each pair $(n_1, n_2)$, we randomly generate 5 instances of problem~\eqref{numerical_example} by sampling the entries of $A$, $B$, and $u$ independently from the standard normal distribution $\mathcal{N}(0,1)$. Note that problem~\eqref{numerical_example} is a special case of problem \eqref{intro:problem} with $f(x,y) = -\| (y + A x) \odot (y + B x) \|^2 + 0.01\|x - u\|^2$, $p(x) = 0.01\|x\|_1 + \delta_{\cB(0,2)}(x)$, $q(y) = 0.1\|y\|_1 + \delta_{[-2, 2]^{n_2}}(y)$, and $c = (c_1, \ldots, c_m)$ defined in \eqref{numerical_constraints}.

We now apply Algorithm~\ref{algo_implement:whole} to solve problem~\eqref{numerical_example} on the randomly generated instances described above. Notice that problem~\eqref{numerical_example} is similar to the one studied in \cite[Section 5]{lu2025first}, except that the constraint $c(y) \leq 0$ is imposed on the inner maximization problem. Consequently, the Lipschitz constant $L_f$ of $f(\cdot,y)$ and the Lipschitz constant $L_{\nabla f}$ of $\nabla f$ over $\cB(0,2) \times [-2, 2]^{n_2}$ are computed as in \cite[Section 5]{lu2025first}. In addition, one can verify that $c_j$ is $L_{c_j}$-smooth over $[-2, 2]^{n_2}$ with $L_{c_j} = e^2$ for all $j \in \{1,...,m\}$. The remaining input parameters for Algorithm~\ref{algo_implement:whole} are set as $C = 0.1$, $\theta = 0.5$, $\gamma = 0.01$, $\sigma = 0.1$, $\underline{L} = 1$, $\rho = 1.25$, $\beta = 10$, $\epsilon = 10^{-2}$.\footnote{
In our numerical experiments, we tested several values of $\theta \in [1/2,1)$ and observed that $\theta = 1/2$ yields the largest reduction in the objective value. The remaining parameters were chosen empirically and perform reasonably well across the tested instances.}
The algorithm is initialized at $(x^0, y^0) = (0, 0)$. Note that for this choice of $(x^0, y^0)$, $y^0$ is the maximizer of the problem $\max_{c(y) \leq 0} \{f(x^0, y) - 0.1\|y\|_1 - \delta_{[-2, 2]^{n_2}}(y)\}$, making it a suitable starting point for $y$. We run the algorithm for 2{,}500 iterations and return the final output denoted by $(x_\epsilon, y_\epsilon)$. Here, $x_\epsilon$ serves as an approximate solution to the outer minimization problem of~\eqref{numerical_example}, while $y_\epsilon$ is an approximate solution to the inner maximization problem $\max_{c(y) \leq 0} \{f(x_\epsilon, y) - 0.1\|y\|_1 - \delta_{[-2, 2]^{n_2}}(y)\}$.

To evaluate the performance of Algorithm~\ref{algo_implement:whole}, we compute the actual final objective value
\[
\Psi(x_\epsilon) = \max_{c(y) \leq 0} \{f(x_\epsilon, y) - 0.1\|y\|_1 - \delta_{[-2, 2]^{n_2}}(y)\} + 0.01 \|x_\epsilon\|_1.
\]
Thanks to the block-separable structure of the problem, this maximization problem can be decomposed into $m$ independent subproblems, each involving a single component of $c$ and $10$ components of $y$.
These subproblems are solved using the \texttt{MATLAB} subroutine \texttt{GlobalSearch}, which is a solver for finding global optima of nonconvex problems. In addition, we compute an approximate final objective value by
\[
\widehat\Psi(x_\epsilon) = f(x_\epsilon, y_\epsilon) - 0.1\|y_\epsilon\|_1 + 0.01\|x_\epsilon\|_1,
\]
using the approximate inner solution $y_\epsilon$ returned by the algorithm.

\begin{figure}[ht]
    \centering
    \includegraphics[width=0.85\linewidth]{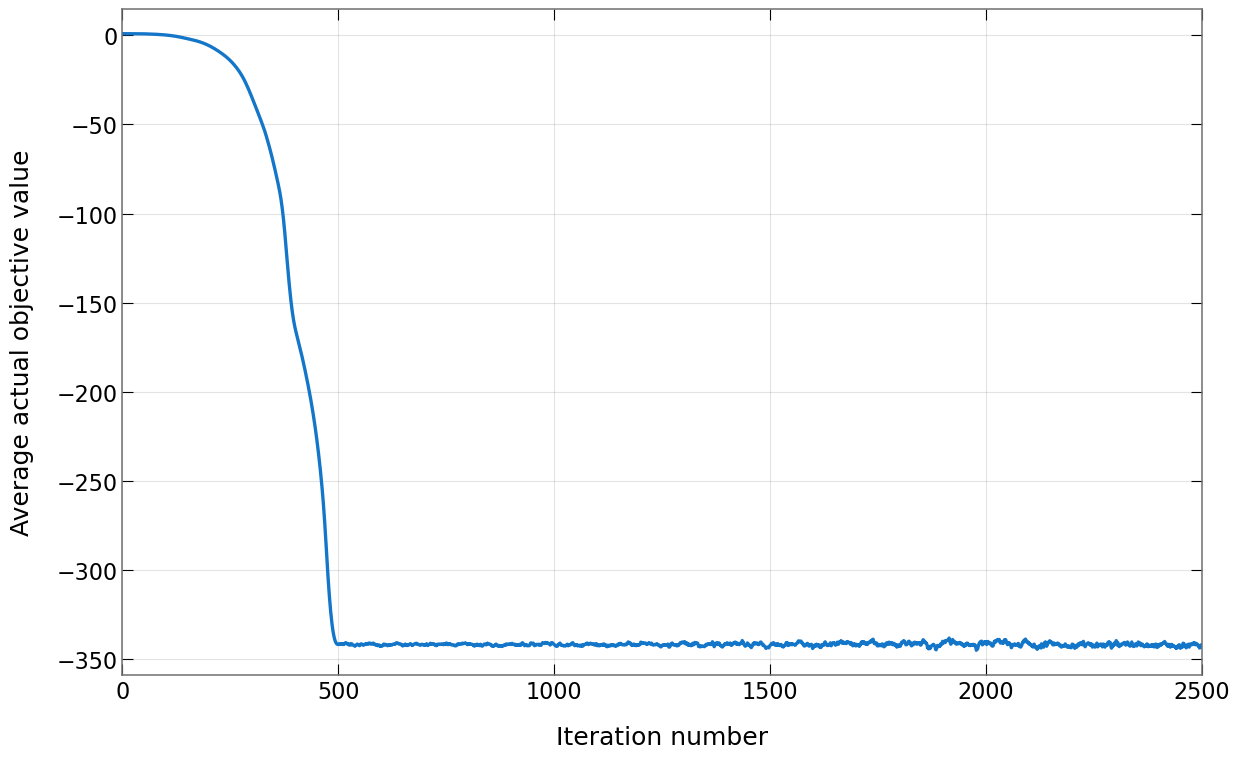}
    \caption{Performance of Algorithm~\ref{algo_implement:whole} with $(n_1,n_2)=(100,100)$}
    \label{fig:obj_descent}
\end{figure}

The computational results on the random instances are presented in Table~\ref{tab:results}. The first two columns list the values of $n_1$ and $n_2$. For each pair $(n_1, n_2)$, the average initial, actual final, and approximate final objective values over five random instances are reported in the remaining columns. From the results, we observe that the approximate solution $x_\epsilon$ significantly reduces the objective value compared to the initial point $x^0$, and that $y_\epsilon$ is a good approximate solution to the inner maximization problem $\max_{c(y) \leq 0} \{f(x_\epsilon, y) - 0.1\|y\|_1 - \delta_{[-2, 2]^{n_2}}(y)\}$.
{In addition, for the five random instances with $(n_1,n_2)=(100,100)$, we plot the average actual objective value in Figure~\ref{fig:obj_descent} to illustrate the performance of Algorithm~\ref{algo_implement:whole}. As observed, the average objective value decreases rapidly in the early stage and then stabilizes, which illustrates the convergence behavior of the proposed method.}

\begin{table}[ht]
\centering
\caption{Numerical results for Algorithm~\ref{algo_implement:whole}}
\label{tab:results}
\begin{tabular}{c c c c c}
\toprule
$n_1$ & $n_2$ & Initial objective value & Actual final value & Approximate final value\\
\midrule
50 & 50 & 0.49 & -164.92 & -165.09 \\
60 & 60 & 0.61 & -279.51 & -279.53 \\
70 & 70 & 0.78 & -241.86 & -241.91 \\
80 & 80 & 0.75 & -298.59 & -298.76 \\
90 & 90 & 0.81 & -260.97 & -261.38 \\
100 & 100 & 1.02 & -341.73 & -341.88 \\
110 & 110 & 1.14 & -404.36 & -404.43 \\
120 & 120 & 1.16 & -482.48 & -482.50 \\
130 & 130 & 1.17 & -467.62 & -467.75 \\
140 & 140 & 1.39 & -590.08 & -590.23 \\
150 & 150 & 1.42 & -579.71 & -579.88 \\
\bottomrule
\end{tabular}
\end{table}

\section{Proof of the main results}\label{sec:proof}
In this section we provide a proof of our main results presented in Sections \ref{sec:localKL_Holder},  \ref{sec:subsolver_SCP},  and \ref{sec:FOD}, which are particularly Lemma~\ref{lambda-bound} and Theorems~\ref{thm:Phi_localHolder},  \ref{thm:Phi_QG}, \ref{thm:sub_well}, \ref{thm:sub_iters}, \ref{thm:outer_bound}, and \ref{thm:operations}. 

\subsection{Proof of the main results in Section~\ref{sec:localKL_Holder}}\label{sec:proof1}

In this subsection, we prove Theorems~\ref{thm:Phi_localHolder} and \ref{thm:Phi_QG}. To this end, we first present several technical lemmas. The following lemma establishes the equivalence between $\max_{y} F(x, y)$ and $\max_{y, z} \bF(x, y, z)$ for any $x \in \cX$, and between $\min_x \Psi(x)$ and $\min_x \bPsi(x)$.

\begin{lemma} \label{lem:equiv_Phi}
Let $F^*, Y^*,\Psi, \bF^*,\overline{Y}^*,\bPsi$ be defined in \eqref{intro_F_Ystar}, \eqref{intro_Phi_star}, and \eqref{intro_PhiOver}. Suppose that Assumption~\ref{ass:Lip_Smo_KL} holds. Then for any {$x \in \cX$}, the following statements hold.
\begin{enumerate} [label=(\roman*)]
\item If $y^* \in Y^*(x)$, then $(y^*, y^*) \in \overline{Y}^*(x)$.
\item $F^*(x) = \bF^*(x)$ and $\Psi(x)=\bPsi(x)$.
\item If $(y^*, z^*) \in \overline{Y}^*(x)$, then $y^* \in Y^*(x)$.
\end{enumerate}
\end{lemma}

\begin{proof}
Fix any $x \in \cX$. For notational convenience, let $\tilde{F}(x, y) = f(x, y) - q(y)$. It then follows from the definitions of $F$ and $\bF$ in \eqref{intro_F_Ystar} and \eqref{intro_potential} that 
\beq \label{lem1_Ftilde}
F(x, y) = \tilde{F}(x, y) - \delta_{c(\cdot) \leq 0}(y), \qquad \bF(x, y, z) = \tilde{F}(x, y) - \delta_{\bc(\cdot, \cdot) \leq 0}(y, z) - \delta_{\cY}(z).
\eeq

We first prove statement (i). Fix any $y^* \in Y^*(x)$. Clearly, $y^* \in \cY$, $c(y^*) \leq 0$, and
\beq \label{lem1_i}
\tilde{F}(x, y) \leq \tilde{F}(x, y^*) \quad \forall y \in \cY \text{ with } c(y) \leq 0.
\eeq
Morever,  by \eqref{intro_cOver} and $c(y^*) \leq 0$, we observe that $\bc(y^*, y^*) = c(y^*) \leq 0$.
Let $(y', z') \in \bR^{n_2} \times  \bR^{n_2}$ be arbitrarily chosen. We claim that $\bF(x, y', z') \leq \bF(x, y^*, y^*)$. Indeed, if $\bc(y', z') > 0$ or $(y', z') \notin \cY \times \cY$, then $\bF(x, y', z') =-\infty$, and the claim holds trivially. Now suppose that $\bc(y', z') \leq 0$ and $(y', z') \in \cY \times \cY$. By these,  \eqref{intro_cOver}, and the $L_{c_i}$-Lipschitz smoothness of each component $c_i$ over $\cY$ (see Assumption~\ref{ass:Lip_Smo_KL}), we deduce that $c(y') \leq \bc(y', z') \leq 0$. This along with $y' \in \cY$ and \eqref{lem1_i} yields $\tilde{F}(x, y') \leq \tilde{F}(x, y^*)$. It then follows from \eqref{lem1_Ftilde}, $\bc(y', z') \leq 0$, $\bc(y^*, y^*) \leq 0$, and $y^*, z' \in \cY$ that
\[
\bF(x, y', z') = \tilde{F}(x, y') \leq \tilde{F}(x, y^*) = \bF(x, y^*, y^*),
\]
and the above claim again holds. By $\bF(x, y', z') \leq \bF(x, y^*, y^*)$ and the arbitrariness of $(y', z')$, we conclude that $(y^*, y^*) \in \overline{Y}^*(x)$. Hence, statement (i) holds.

We next prove statement (ii). By Assumption~\ref{ass:Lip_Smo_KL}, there exists at least one $\hat{y}^* \in Y^*(x)$, and moreover, $F(x, \hat{y}^*)$ is finite. Using this and statement (i), we see that $(\hat{y}^*, \hat{y}^*) \in \overline{Y}^*(x)$. Notice that $\hat{y}^* \in \cY$ and $\bc(\hat{y}^*, \hat{y}^*)=c(\hat{y}^*) \leq 0$. By these, \eqref{intro_F_Ystar}, \eqref{intro_PhiOver}, and \eqref{lem1_Ftilde}, we obtain that 
\[
F^*(x) = F(x, \hat{y}^*) = \tilde{F}(x, \hat{y}^*) = \bF(x, \hat{y}^*, \hat{y}^*) = \bF^*(x).
\]
It follows from this, \eqref{intro_Phi_star}, and \eqref{intro_PhiOver} that $\Psi(x)=\bPsi(x)$ holds. 
This proves statement (ii).

We finally prove statement (iii). Fix any $(y^*, z^*) \in \overline{Y}^*(x)$. By this, \eqref{intro_potential}, \eqref{intro_PhiOver}, and \eqref{lem1_Ftilde}, we observe that $(y^*, z^*) \in \cY \times \cY$, $\bc(y^*, z^*) \leq 0$, and $\bF^*(x)=\bF(x, y^*, z^*)=\tilde{F}(x, y^*)$. Using these relations,  \eqref{intro_cOver}, and the $L_{c_i}$-Lipschitz smoothness of each component $c_i$ over $\cY$, we deduce that $c(y^*) \leq \bc(y^*, z^*) \leq 0$. By this, 
 \eqref{lem1_Ftilde}, and $y^* \in \cY$, one has $F(x,y^*)=\tilde{F}(x, y^*)$, which together with $\bF^*(x)=\tilde{F}(x, y^*)$ and $F^*(x) = \bF^*(x)$ (see statement (ii)) implies that $F(x,y^*)=F^*(x)$. 
Hence, $y^*\in Y^*(x)$ and statement (iii) holds. 
\end{proof}

\begin{remark} \label{minimax-equiv} 
In view of \eqref{intro_F_Ystar}, \eqref{intro_PhiOver}, and  Lemma \ref{lem:equiv_Phi}(ii), we observe that 
\[
\max_{y} F(x, y) = \max_{y, z} \bF(x, y, z) \qquad \forall x \in \cX.
\]
Consequently, when interpreted as minimization problems, $\min_x \{\max_{y} F(x, y)+p(x)\}$ and \\ $\min_x \{ \max_{y, z} \bF(x, y, z)+p(x)\}$ have identical objective functions and are thus equivalent. Moreover, from \eqref{intro_F_Ystar}, \eqref{intro_potential}, and \eqref{intro_PhiOver}, these problems correspond to \eqref{intro:problem} and \eqref{equiv-prob}, respectively. Therefore, the original minimax problem \eqref{intro:problem} is equivalent to the lifted minimax problem \eqref{equiv-prob}.
\end{remark}

The next lemma presents properties of the limiting slope and error bound of a proper closed function.  Its proof can be found in \cite[Proposition 4.6]{drusvyatskiy2015curves} and \cite[Lemma 2.5]{drusvyatskiy2015curves}.

\begin{lemma} \label{lem:errbnd}
Let $\phi: \bR^n \to \overline{\bR}$ be a proper closed function, and a point $\bu\in\dom\,\phi$ be given. Then the following statements hold.
\begin{itemize}
\item[(i)] $\overline{|\nabla \phi|}(\bu)=\dist(0,\partial \phi(\bu))$. 
\item[(ii)] Suppose there exist constants $\alpha < \phi(\bu)$ and $r, K > 0$ such that
\[
\alpha < \phi(u) \leq \phi(\bu) \   \text{ and } \   \|u-\bu\| \leq K  \quad \Longrightarrow \quad  |\nabla \phi |(u)  \geq r. 
\]
If, in addition, $\phi(\bu)-\alpha < Kr$, then 
\[
\dist(\bu, \cS_\alpha)\leq r^{-1}(\phi(\bu)-\alpha), \ \   \text{ where} \ \  \cS_\alpha:= \{u: \phi(u) \leq \alpha\}.
\]
\end{itemize}
\end{lemma}

The lemma below establishes a local $(1-\theta)^{-1}$-growth property of $\bF(x, \cdot, \cdot)$ for any $x \in \cX$, which can also be obtained from \cite[Theorem 3.7]{drusvyatskiy2021nonsmooth}. For completeness, we provide a self-contained proof with minimal reliance on the literature.

\begin{lemma} \label{lem:QG_pot}
Suppose that Assumption~\ref{ass:Lip_Smo_KL} holds. Then it holds that for any {$x \in \cX$},  
\beq \label{lem2_QG_pot}
\bF^*(x) - \bF(x, y, z) \geq (C(1-\theta))^{\frac{1}{1-\theta}} \dist((y, z), \overline{Y}^*(x))^{\frac{1}{1-\theta}} \qquad \forall (y, z) \in \overline{\cL}(x).
\eeq
\end{lemma}

\begin{proof}
Fix any $x \in \cX$ and $(\by, \bz) \in \overline{\cL}(x)$. If $(\by, \bz) \notin \dom\,\bF(x, \cdot, \cdot)$, then $\bF(x, \by,\bz)=-\infty$, and hence relation \eqref{lem2_QG_pot} holds trivially at $(y,z)=(\by, \bz)$.  We now suppose that $(\by, \bz) \in \dom\,\bF(x, \cdot, \cdot)$. For notational convenience, let 
\begin{align}
& \phi_x(y, z) = -\bF(x, y, z), \quad \phi_x^* = -\bF^*(x),  \label{phix}\\
& g(y, z) = (\phi_x(y, z) - \phi_x^*)^{1-\theta}, \quad  K_x = 1+(C(1-\theta))^{-1}(\gamma \dist(0, \partial \bPsi(x))^\sigma)^{1-\theta}. \label{g-fun}
\end{align}
Then, one can see from \eqref{poten_KL}, \eqref{poten_Lev}, and $(\by, \bz) \in \overline{\cL}(x)$ that
\beq \label{lem2_phi_KL}
0<C(\phi_x(\by, \bz) - \phi_x^*)^\theta \leq \dist(0, \partial \phi_x(\by, \bz)).
\eeq
By the definitions of slope and limiting slope in \eqref{slope} and \eqref{limit-slope}, one can observe that $|\nabla \phi_x|(\by, \bz) \geq \overline{|\nabla \phi_x|}(\by, \bz)$. Also, it follows from Lemma \ref{lem:errbnd}(i) that $\overline{|\nabla \phi_x|}(\by, \bz)=\dist(0, \partial \phi_x(\by, \bz))$. Hence, we obtain that 
\beq \label{slope-lowbnd}
|\nabla \phi_x|(\by, \bz) \geq \dist(0, \partial \phi_x(\by, \bz)).
\eeq
 In addition, by \eqref{g-fun} and the concavity of the function $t^{1-\theta}$ in $[0,\infty)$ due to $\theta \in [1/2,1)$, one has 
\begin{align*}
g(y, z) &\overset{\eqref{g-fun}}{=} (\phi_x(y, z) - \phi_x^*)^{1-\theta} \leq (\phi_x(\by, \bz) - \phi_x^*)^{1-\theta} + (1-\theta) (\phi_x(\by, \bz) - \phi_x^*)^{-\theta} (\phi_x(y, z)-\phi_x(\by, \bz)) \\
&=  g(\by, \bz) + (1-\theta) (\phi_x(\by, \bz) - \phi_x^*)^{-\theta} (\phi_x(y, z)-\phi_x(\by, \bz)).
\end{align*}
Using this,  $\theta \in [1/2,1)$, $\phi_x(\by, \bz) > \phi_x^*$,  and the definition of slope in  \eqref{slope}, we obtain that 
\begin{align*}
|\nabla g|(\by, \bz) &\overset{ \eqref{slope}}{=}  \limsup_{(y,z)\to (\by,\bz)} \frac{\big(g(\by,\bz)-g(y,z)\big)_+}{\|(\by,\bz)-(y,z)\|} \\
& \geq \limsup_{(y,z)\to (\by,\bz)}  \frac{\big((1-\theta) (\phi_x(\by, \bz) - \phi_x^*)^{-\theta} (\phi_x(\by, \bz)-\phi_x(y, z))\big)_+}{\|(\by,\bz)-(y,z)\|} \\
& =(1-\theta) (\phi_x(\by, \bz) - \phi_x^*)^{-\theta}   \limsup_{(y,z)\to (\by,\bz)} \frac{\big(\phi_x(\by, \bz)-\phi_x(y, z)\big)_+}{\|(\by,\bz)-(y,z)\|} \\
&\overset{ \eqref{slope}}{=}(1-\theta) (\phi_x(\by, \bz) - \phi_x^*)^{-\theta} |\nabla \phi_x|(\by, \bz) \\ 
& \overset{\eqref{slope-lowbnd}}{\geq} (1-\theta) (\phi_x(\by, \bz) - \phi_x^*)^{-\theta}  \dist(0, \partial \phi_x(\by, \bz)) \overset{\eqref{lem2_phi_KL}}{\geq} C(1-\theta).
\end{align*}
By this relation and the arbitrariness of $(\by, \bz) \in \overline{\cL}(x) \cap \dom\,\bF(x, \cdot, \cdot)$, we conclude that $|\nabla g|(y, z) \geq  C(1-\theta)$ holds for all $(y, z) \in \overline{\cL}(x) \cap \dom\,\bF(x, \cdot, \cdot)$. In addition, by the definitions of $g$ and $\overline{\cL}(x)$ along with the fact $(\by, \bz) \in \overline{\cL}(x) \cap \dom\,\bF(x, \cdot, \cdot)$, one can observe that for any $(y, z)$ with $0 < g(y, z) \leq g(\by, \bz)$, we have $(y, z) \in \overline{\cL}(x) \cap \dom\,\bF(x, \cdot, \cdot)$, and hence $|\nabla g|(y, z) \geq C(1-\theta)$. Also, notice from \eqref{poten_Lev} and $(\by, \bz) \in \overline{\cL}(x)$ that 
\[
0 < \bF^*(x) - \bF(x, \by, \bz) \leq \gamma\,\dist(0,\partial\bPsi(x))^\sigma,
\]
which along with \eqref{phix} and \eqref{g-fun} implies that  
\begin{align*}
0< g(\by, \bz) &\overset{\eqref{g-fun}}{=} (\phi_x(\by, \bz) - \phi_x^*)^{1-\theta} \overset{\eqref{phix}}{=} (\bF^*(x) - \bF(x, \by, \bz))^{1-\theta} \\ 
& \leq \big(\gamma\,\dist(0,\partial\bPsi(x))^\sigma\big)^{1-\theta} \overset{\eqref{g-fun}}{<} K_x C(1-\theta).
\end{align*}
 In view of these, it follows from  \eqref{phix}, \eqref{g-fun}, and Lemma \ref{lem:errbnd} with $\phi=g$, $\bu=(\by,\bz)$, $\alpha=0$, $r=C(1-\theta)$, $K=K_x$, and $\cS_\alpha=\overline{Y}^*(x)$ that 
\[
\dist((\by, \bz), \overline{Y}^*(x)) \overset{\text{Lemma} \ \ref{lem:errbnd}}\leq \frac{1}{C(1-\theta)} g(\by, \bz) \overset{\eqref{phix}\eqref{g-fun}}{=} \frac{1}{C(1-\theta)} (\bF^*(x) - \bF(x, \by, \bz))^{1-\theta},
\]
which implies that relation \eqref{lem2_QG_pot} holds at $(y,z)=(\by, \bz)$. By this and the arbitrariness of $(\by, \bz) \in \overline{\cL}(x) \cap \dom\,\bF(x, \cdot, \cdot)$, one can conclude that relation \eqref{lem2_QG_pot} holds for all $(y, z) \in \overline{\cL}(x) \cap \dom\,\bF(x, \cdot, \cdot)$. This completes the proof. 
\end{proof}

The following lemma provides a relationship between $\dist\big((y, z), \overline{Y}^*(x)\big)$ and $\dist\big(0, \partial_{(y, z)} \bF(x, y, z)\big)$, following directly from \eqref{poten_KL} and \eqref{lem2_QG_pot}. 

\begin{lemma} \label{lem:EB_pot}
Suppose that Assumption~\ref{ass:Lip_Smo_KL} holds. Then it holds that for any {$x \in \cX$}, 
\[
\dist((y, z), \overline{Y}^*(x)) \leq (1-\theta)^{-1} {C^{-\frac{1}{\theta}}} \dist\big(0, \partial_{(y, z)} \bF(x, y, z)\big)^{\frac{1-\theta}{\theta}} \qquad \forall (y, z) \in \overline{\cL}(x).
\]
\end{lemma}

{The lemma below establishes the local {generalized} H\"{o}lder smoothness of $\bF^*$. It is analogous to \cite[Theorem 1]{lu2025first}, whose proof is based on \cite[Lemma~4]{lu2025first}.  Note that \cite[Lemma~4]{lu2025first} is analogous to Lemma \ref{lem:EB_pot}. Therefore, the proof of this lemma follows directly from Lemma~\ref{lem:EB_pot}, together with arguments similar to those used in the proof of \cite[Theorem~1]{lu2025first}, and is thus omitted.}

\begin{lemma} \label{lem:PhiOver_localHolder}
Let $\epsilon>0$ be given and $\overline{\cU}_\epsilon = \{x \in \cX: \dist(0, \partial \bPsi(x)) > \epsilon\}$. Suppose that Assumption~\ref{ass:Lip_Smo_KL} holds. Then the following statements hold.
\begin{itemize}
\item[(i)] $\partial^\rC_\cX \bF^*(x)$ is a singleton for all $x \in \overline{\cU}_\epsilon$, and $\bF^*$ is differentiable on $\overline{\cU}_\epsilon \cap \operatorname{int}(\cX)$. 
\item[(ii)]  For any $x, x' \in \overline{\cU}_\epsilon \cap \operatorname{int}(\cX)$ satisfying $\|x - x'\| \leq \gamma \epsilon^\sigma/(2L_f)$, we have
\[
\|\nabla \bF^*(x) - \nabla \bF^*(x')\| \leq L_{\nabla f}\|x - x'\| + (1-\theta)^{-1} {C^{-1/\theta}}L_{\nabla f}^{1/\theta}\, \|x- x'\|^{\frac{1-\theta}{\theta}}.    
\]
\item[(iii)] For any $x, x' \in \overline{\cU}_\epsilon$ satisfying $\|x - x'\| \leq \gamma \epsilon^\sigma/(4L_f)$, we have
\[
\|\nabla^\rC_\cX \bF^*(x) - \nabla^\rC_\cX \bF^*(x')\| \leq L_{\nabla f}\|x - x'\| + (1-\theta)^{-1} {C^{-1/\theta}}L_{\nabla f}^{1/\theta}\, \|x- x'\|^{\frac{1-\theta}{\theta}}.    
\]
\item[(iv)] It holds that
\[
\nabla^\rC_\cX \bF^*(x) = \nabla_x f(x, y^*) \qquad \forall x \in \overline{\cU}_\epsilon, \ y^* \in Y^*(x).
\]
\end{itemize}
\end{lemma}

We are now ready to prove Theorems~\ref{thm:Phi_localHolder} and \ref{thm:Phi_QG}.

\begin{proof}[\textbf{Proof of Theorem~\ref{thm:Phi_localHolder}}] 
The conclusion of this theorem directly follows from Lemmas~\ref{lem:equiv_Phi} and \ref{lem:PhiOver_localHolder}.
\end{proof}

\begin{proof}[\textbf{Proof of Theorem~\ref{thm:Phi_QG}}] 
Fix any $x \in \cX$ and $y \in \cL(x)$. If $y \notin \dom\,F(x, \cdot)$,  then $F(x, y) =-\infty$, and hence the conclusion holds trivially. We now suppose $y \in \dom\,F(x, \cdot)$. It follows from this and the definitions of $F$ and $\bc$ in \eqref{intro_F_Ystar} and \eqref{intro_cOver} that $y \in \cY$ and $\bc(y, y) = c(y) \leq 0$, which implies that $(x, y, y) \in \dom\,\bF(x, \cdot, \cdot)$ and $F(x, y) = \bF(x, y, y)$. Recall from Lemma~\ref{lem:equiv_Phi} that $F^*(x) = \bF^*(x)$. By these,  \eqref{poten_Lev},  \eqref{Lev_set}, and $y \in \cL(x)$,  one has $(y, y) \in \overline{\cL}(x)$. Using these relations and Lemma~\ref{lem:QG_pot}, we have
\beq \label{lem5_QGtmp}
F^*(x) - F(x, y) = \bF^*(x) - \bF(x, y, y) \geq (C(1-\theta))^{\frac{1}{1-\theta}} \dist((y, y), \overline{Y}^*(x))^{\frac{1}{1-\theta}}.
\eeq

We next show that $\dist(y, Y^*(x)) \leq \dist((y, y), \overline{Y}^*(x))$. Notice from Assumption~\ref{ass:Lip_Smo_KL} and Lemma~\ref{lem:equiv_Phi} that $\overline{Y}^*(x)$ is a nonempty closed set. Hence, there exists $(y^*, z^*) \in \overline{Y}^*(x)$ such that $\|(y, y) - (y^*, z^*)\| = \dist((y, y), \overline{Y}^*(x))$. Since $(y^*, z^*) \in \overline{Y}^*(x)$, it follows from Lemma~\ref{lem:equiv_Phi} that $y^* \in Y^*(x)$, which implies that $\dist(y, Y^*(x)) \leq \|y - y^*\|$. In view of these, one has
\[
\dist(y, Y^*(x)) \leq \|y - y^*\| \leq \|(y, y) - (y^*, z^*)\| = \dist((y, y), \overline{Y}^*(x)),
\]
and hence $\dist(y, Y^*(x)) \leq \dist((y, y), \overline{Y}^*(x))$ holds as desired. The conclusion \eqref{lem5_Phi_QG} directly follows from this and \eqref{lem5_QGtmp}.
\end{proof}

\subsection{Proof of the main results in Section~\ref{sec:subsolver_SCP}}\label{sec:proof2}

In this subsection we prove Theorems~\ref{thm:sub_well} and \ref{thm:sub_iters}. 

\begin{proof}[\textbf{Proof of Theorem~\ref{thm:sub_well}}] 
We first prove statement (i) of Theorem~\ref{thm:sub_well} by induction. Suppose that a point $z^k \in \dom\,q$ satisfying $c(z^k) \leq 0$ is already generated in Algorithm~\ref{algo_implement:subsolver} for some $k\geq 0$. Let us fix any $i \geq 0$, 
and define
\beq \label{Q-fun}
Q_{k,i}(z)=\langle \nabla g(z^k), z \rangle + \frac{L_{k, i}}{2}\|z - z^k\|^2 + q(z) +\delta_{\bc(\cdot,z^k)\leq 0}(z).
\eeq
Since $z^k \in \dom\,q$ and $\bc(z^k, z^k) = c(z^k) \leq 0$, one has $z^k\in\dom\,Q_{k,i}$. By this and the assumption that $q$ is a closed convex function, one can observe that $Q_{k,i}$ is a proper closed strongly convex function. It follows that the problem $\min_z Q_{k,i}(z)$ has a unique optimal solution. Note that subproblem \eqref{surro_argmin} is equivalent to $z^{k+1,i}=\mathop{\arg\min}_z Q_{k,i}(z)$. Hence,   \eqref{surro_argmin}  has a unique optimal solution $z^{k+1, i}\in\dom\,q$ satisfying $\bc(z^{k+1,i}, z^k) \leq 0$. By these, $z^k \in \dom\,q$, and the Lipschitz smoothness of $c$ on $\dom\,q$, one can conclude that $c(z^{k+1, i}) \leq \bc(z^{k+1,i}, z^k) \leq 0$. 
{In addition, since $z^k \in \dom\,q$ and $c(z^k) \leq 0$, it follows that $z^k\in\cS$, where $\cS$ is defined in \eqref{intro_domain}. By this fact and Theorem~\ref{thm:UniSlater}, the Slater's condition holds for \eqref{surro_argmin}, that is, there exists a point $\tilde{z} \in \dom\,q$ such that $\bc(\tilde{z}, z^k) < 0$.}
Hence, it follows from \cite[Corollary 28.2.1, Theorem 28.3]{rockafellar1997convex} that subproblem \eqref{surro_argmin} has an optimal Lagrange multiplier $\lambda^{k, i}$. We next show that the inner loop terminates in at most $\bar{i}+1$ iterations and outputs a point $z^{k+1} \in \dom\,q$ with $c(z^{k+1}) \leq 0$. To this end, suppose for contradiction that the inner loop runs for more than $\bar{i}+1$ iterations. Then one can observe from Algorithm~\ref{algo_implement:subsolver} that
\beq \label{thm2_contra}
h(z^{k+1,\bar{i}}) > h(z^k)- \frac{\beta}{2} \|z^{k+1,\bar{i}} - z^k\|^2. 
\eeq 
Since $z^{k+1,\bi}=\mathop{\arg\min}_z Q_{k,\bi}(z)$ and $Q_{k,\bi}$ is strongly convex with modulus $L_{k, \bi}$, it follows that $Q_{k,\bi}(z^{k+1,\bar{i}}) \leq Q_{k,\bi}(z^{k})-L_{k, \bi}\|z^{k+1,\bar{i}}-z^k\|^2/2$, which together with $z^k, z^{k+1, \bi}\in\dom\,Q_{k,\bi}$ and the definition of $Q_{k,\bi}$ in \eqref{Q-fun} yields
\[
\langle \nabla g(z^k), z^{k+1, \bi} \rangle + \frac{L_{k, \bi}}{2} \|z^{k+1, \bi} - z^k\|^2 + q(z^{k+1, \bi})  \leq \langle \nabla g(z^k), z^{k} \rangle + q(z^{k}) - \frac{L_{k, \bi}}{2}\|z^{k+1, \bi} - z^k\|^2.
\]
By this,  \eqref{pr:KL_opt}, and the $L$-smoothness of $g$,  one has
\begin{align}
h(z^{k+1, \bi}) &\overset{\eqref{pr:KL_opt}}{=} g(z^{k+1, \bi}) + q(z^{k+1, \bi})
\leq g(z^k) + \langle \nabla g(z^k), z^{k+1, \bi} - z^k \rangle + \frac{L}{2} \|z^{k+1, \bi} - z^k\|^2 + q(z^{k+1, \bi}) \nn \\
&\leq g(z^k) + q(z^k) - \Big(L_{k,\bi} - \frac{L}{2}\Big) \|z^{k+1, \bi} - z^k\|^2  \overset{\eqref{pr:KL_opt}}{=} h(z^k) - \Big(L_{k, \bi} - \frac{L}{2}\Big) \|z^{k+1,\bi} - z^k\|^2. \label{thm2_innerBound}
\end{align}
Notice from $L_{k,\bar{i}} = \underline{L}\rho^{\bar{i}}$ and the definition of $\bar{i}$ that $L_{k,\bar{i}} \geq (\beta + L)/2$. This and \eqref{thm2_innerBound} lead to $h(z^{k+1,\bar{i}})  \leq h(z^k)-\beta \|z^{k+1,\bar{i}} - z^k\|^2/2$, which contradicts \eqref{thm2_contra}. Hence, the inner loop terminates in at most $\bar{i}+1$ iterations. Moreover, it outputs a point $z^{k+1} \in \dom\,q$ satisfying $c(z^{k+1}) \leq 0$ due to $z^{k+1, i} \in \dom\,q$ and $c(z^{k+1, i}) \leq 0$ for each $i$. This together with the fact that $z^0 \in \dom\,q$ and $c(z^0) \leq 0$ implies that the induction is complete. Hence, statement (i) holds as desired.

We next  prove statement (ii) of Theorem~\ref{thm:sub_well}.  By the definition of $L_k$ and statement (i), one can see that $L_k = L_{k, i}$ for some $0\leq i \leq \bar{i}$, which together with the definition of $\bar{i}$ implies that \eqref{Lk-bound} holds. In addition, notice from Algorithm~\ref{algo_implement:subsolver} that $(z^{k+1},\lambda^k)$ is a pair of optimal solution and Lagrange multiplier of the problem 
 \[
\min_z \Big\{ \langle \nabla g(z^k), z \rangle + \frac{L_k}{2}\|z - z^k\|^2 + q(z): \bc(z, z^k) \leq 0\Big\}.
\]
The relations  \eqref{thm2_KKT_CompSlack} and \eqref{thm2_KKT_Stationary} then follow from the KKT conditions of this problem at $(z^{k+1},\lambda^k)$.
\end{proof}

We now turn to the proof of Lemma~\ref{lambda-bound}.

\begin{proof}[\textbf{Proof of Lemma~\ref{lambda-bound}}] 
Let $(z^{k+1},\lambda^k,L_k)$ be generated in the $k$th outer iteration of Algorithm~\ref{algo_implement:subsolver} for some $k\geq 0$. Observe from Algorithm~\ref{algo_implement:subsolver}  that
\[
\begin{aligned}
z^{k+1} = &\mathop{\arg\min}_z \Big\{ \langle \nabla g(z^k), z \rangle + \frac{L_k}{2}\|z - z^k\|^2 + q(z) \Big\} \\
&\qquad\text{s.t.}  \quad \bc(z, z^k) \leq 0,
\end{aligned}
\]
and $\lambda^k$ is its associated optimal Lagrange multiplier.  It follows that 
\beq \label{KKT-cond}
\lambda^k \geq 0, \quad \langle \lambda^k, \bc(z^{k+1}, z^k) \rangle=0, \quad
z^{k+1} = \mathop{\arg\min}_z \widetilde{L}(z, \lambda^k), 
\eeq
where 
\[
\widetilde{L}(z, \lambda) =\phi(z) + \langle \lambda, \bc(z, z^k) \rangle \ \  \text{with} \ \ \phi(z)= \langle \nabla g(z^k), z \rangle + \frac{L_k}{2}\|z - z^k\|^2 + q(z). 
\]
Notice from Theorem \ref{thm:sub_well}(i) that $z^k\in \cS$.  By this and Theorem~\ref{thm:UniSlater},  there exists $\hat{y}^k \in \cY$ such that $\bc(\hat{y}^k, z^k) \leq -\zeta < 0$.  Using this and \eqref{KKT-cond},  we have
\[
\phi(z^{k+1}) = \widetilde{L}(z^{k+1}, \lambda^k)= \min_z \widetilde{L}(z, \lambda^k) \leq \widetilde{L}(\hat{y}^k, \lambda^k) = \phi(\hat{y}^k) + \langle \lambda^k, \bc(\hat{y}^k, z^k) \rangle 
\leq \phi(\hat{y}^k) - \zeta \|\lambda^k\|_1.
\]
It follows from this and $L_k \leq \bar{L}$ (see Theorem~\ref{thm:sub_well}(ii)) that
\[
\begin{aligned}
\|\lambda^k\|_1 &\leq \zeta^{-1} \big(\phi(\hat{y}^k) - \phi(z^{k+1})\big) \leq \zeta^{-1} \Big(\langle \nabla g(z^k), \hat{y}^k-z^{k+1} \rangle + \frac{L_k}{2}\|\hat{y}^k-z^{k+1}\|^2 + q(\hat{y}^k) - q(z^{k+1})\Big) \\
&\leq \zeta^{-1} \Big( \|\nabla g(z^k)\| \|\hat{y}^k-z^{k+1}\| + \frac{\bar{L}}{2}\|\hat{y}^k-z^{k+1}\|^2 + q(\hat{y}^k) - q(z^{k+1}) \Big).
\end{aligned}
\]
Using this, $\hat{y}^k, z^k, z^{k+1} \in \cY$, and the definitions of $G$, $D_{\cY}$ and $M_q$, we see that \eqref{prop1_explicitBound} holds.
\end{proof}

In the remainder of this subsection, we prove Theorem~\ref{thm:sub_iters}. To this end, we first establish several technical lemmas.

The next lemma provides a bound on $\dist(0, \partial \bh(z^{k+1}, z^k))$ in terms of $\|z^{k+1} - z^k\|$.

\begin{lemma} \label{lem:distBound}
Suppose that Assumptions \ref{ass:Lip_Smo_KL}, \ref{ass:MFCQ}, and \ref{assum:g-h} hold.  Let $\bh$, $\{L_{c_j}\}_{j=1}^m$, and $L$ be given in \eqref{subpr_overh} and  Assumptions \ref{ass:Lip_Smo_KL} and \ref{assum:g-h}, respectively. Suppose that  $z^{k+1}$, $z^k$, $L_k$, and $ \lambda^k$ are generated by Algorithm~\ref{algo_implement:subsolver} for some $k \geq 0$. Then it holds that
\beq \label{subdiff-bound}
\dist(0, \partial \bh(z^{k+1}, z^k)) \leq \Big( (L+L_k)^2 + 4\|\lambda^k\|_\infty^2 \big(\sum_{j=1}^m L_{c_j}\big)^2 \Big)^{\frac{1}{2}} \|z^{k+1} - z^k\|.
\eeq
\end{lemma}

\begin{proof}
Let us fix any $(z,w) \in \dom\,\bh$. It follows from the definition of $\bh$ in \eqref{subpr_overh} that
\begin{align}
&\partial  \bh(z,w) \supseteq \hat{\partial} \bh(z,w) \supseteq
\begin{pmatrix}
    \nabla g(z) + \hat{\partial} q(z) \\
    \hat{\partial} \delta_{\cY}(w)   
\end{pmatrix}
+ \hat{\partial} \delta_{\,\bc(\cdot, \cdot) \leq 0} (z,w) \nn \\
&=
\begin{pmatrix}
    \nabla g(z) + \partial q(z) \\
    \cN_{\cY}(w)
\end{pmatrix}
+ \widehat{\cN}_{\bc(\cdot, \cdot) \leq 0}(z,w) \nn \\
&\supseteq \left\{
\begin{pmatrix}
\nabla g(z) + \partial q(z) + \sum_{j=1}^m \lambda_j \Big(\nabla c_j(w) + L_{c_j}(z-w)\Big) \\
\cN_{\cY}(w) + \sum_{j=1}^m \lambda_j \Big(\nabla^2 c_j(w) (z-w) - L_{c_j}(z-w) \Big)
\end{pmatrix} : \lambda \in \cN_{-\bR^m_+}(\bc(z,w)) \right\}, \label{lem2_subdiff}
\end{align}
where the second relation follows from \eqref{subpr_overh} and \cite[Exercise 8.8, Proposition 10.5, Corollary 10.9]{rockafellar2009variational}, the third relation uses \cite[Proposition 8.12, Exercise 8.14]{rockafellar2009variational} together with the convexity of $q$ and $\cY$, and the last relation follows from \eqref{intro_cOver}, \cite[Theorem 6.14]{rockafellar2009variational}, and $\widehat{\cN}_{-\bR^m_+}(\cdot) = \cN_{-\bR^m_+}(\cdot)$. 

Observe from \eqref{subpr_overh}, Theorem~\ref{thm:sub_well}, and Algorithm~\ref{algo_implement:subsolver} that $(z^{k+1}, z^k) \in \dom\,\bh$. By this and  \eqref{lem2_subdiff}, one has that for any $\lambda \in \cN_{-\bR^m_+}(\bc(z^{k+1}, z^k))$,
\[
\partial \bh(z^{k+1}, z^k) \supseteq
\begin{pmatrix}
    \nabla g(z^{k+1}) + \partial q(z^{k+1}) + \sum_{j=1}^m \lambda_j\big(\nabla c_j(z^k) + L_{c_j}(z^{k+1} - z^k)\big) \\ 
    \sum_{j=1}^m \lambda_j\big(\nabla^2 c_j(z^k) (z^{k+1}-z^k) - L_{c_j}(z^{k+1}-z^k)\big)
\end{pmatrix}.
\]
Notice from \eqref{intro_cOver} and \eqref{thm2_KKT_CompSlack} that $\lambda^k \in \cN_{-\bR^m_+}(\bc(z^{k+1}, z^k))$. In addition, observe from \eqref{thm2_KKT_Stationary} that
\[
\nabla g(z^{k+1}) - \nabla g(z^k) - L_k(z^{k+1}-z^k) \in \nabla g(z^{k+1}) + \partial q(z^{k+1}) + \sum_{j=1}^m \lambda^k_j\big(\nabla c_j(z^k) + L_{c_j}(z^{k+1}-z^k)\big).
\]
In view of these, one has
\beq \label{lem3_subdiffElem}
\begin{pmatrix}
    \nabla g(z^{k+1}) - \nabla g(z^k) - L_k(z^{k+1}-z^k) \\
    \sum_{j=1}^m \lambda^k_j\big(\nabla^2 c_j(z^k) (z^{k+1}-z^k) - L_{c_j}(z^{k+1}-z^k)\big)
\end{pmatrix}
\in \partial \bh(z^{k+1}, z^k).
\eeq
Notice from the $L$-smoothness of $g$ that
\beq \label{lem3_b1}
\begin{aligned}
\|\nabla g(z^{k+1})\!-\!\nabla g(z^k)\!-\!L_k(z^{k+1}\!-\!z^k)\| \leq \| \nabla g(z^{k+1})\!-\!\nabla g(z^k) \| + \|L_k(z^{k+1}\!-\!z^k)\| \leq (L\!+\!L_k) \|z^{k+1}\!-\!z^k\|.
\end{aligned}
\eeq
On the other hand, since $\lambda^k \in \bR^m_+$, we have
\beq \label{lem3_b2}
\begin{aligned}
&\big\|\sum_{j=1}^m \lambda^k_j\big(\nabla^2 c_j(z^k) (z^{k+1}-z^k) - L_{c_j}(z^{k+1}-z^k)\big)\big\| \leq \sum_{j=1}^m \lambda^k_j \|\nabla^2 c_j(z^k) (z^{k+1}-z^k) - L_{c_j}(z^{k+1}-z^k)\|  \\
&\leq \sum_{j=1}^m\!\lambda^k_j \big(\|\nabla^2 c_j(z^k) (z^{k+1}\!-\!z^k)\| \!+\! \|L_{c_j}(z^{k+1}\!-\!z^k)\|\big) \leq 2\!\sum_{j=1}^m \!\lambda^k_j L_{c_j}\|z^{k+1}\!-\!z^k\| \leq 2 \|\lambda^k\|_\infty \! \sum_{j=1}^m\!L_{c_j} \|z^{k+1}\!-\!z^k\|,
\end{aligned}
\eeq
where the third inequality follows from $z^k, z^{k+1} \in \dom\,q$, the convexity of  $\dom\,q$, and $L_{c_j}$-smoothness of $c_j$ over $\dom\,q$.  Combining \eqref{lem3_subdiffElem}, \eqref{lem3_b1}, and \eqref{lem3_b2} yields \eqref{subdiff-bound}, and hence the conclusion holds.
\end{proof}

The next lemma establishes the convergence rate of Algorithm~\ref{algo_implement:subsolver} under suitable assumptions.

\begin{lemma}
Suppose that Assumptions \ref{ass:Lip_Smo_KL}, \ref{ass:MFCQ}, and \ref{assum:g-h} hold. Let
$\{(z^\ell, \lambda^{\ell-1})\}_{\ell=1}^k$ be generated by Algorithm~\ref{algo_implement:subsolver} for some $k \geq 1$, and let $\theta, \eta, \omega, \alpha, C'$ be given in \eqref{KL_subpr}, \eqref{lem9_omega}, and \eqref{lem4:alpha_Cp}, respectively. Then the following statements hold.
\begin{enumerate} [label=(\roman*)]
    \item If $\theta = 1/2$, then
    \beq \label{lem4:case1}
    h(z^k) - h^* \leq \eta (1+\alpha)^{-k}.
    \eeq
    \item If $\theta \in (1/2, 1)$, then
    \beq \label{lem4:case2}
    h(z^k) - h^* \leq \Big( \frac{1}{C'(2\theta-1)} \Big)^{\frac{1}{2\theta-1}} k^{-\frac{1}{2\theta-1}}.
    \eeq
\end{enumerate}
\end{lemma}

\begin{proof}
For notational convenience, let $r_\ell := h(z^{\ell}) - h^*$ for all $0 \leq \ell \leq k$. If $h(z^k)= h^*$, then relations \eqref{lem4:case1} and \eqref{lem4:case2} clearly hold. For the remainder of the proof, suppose that $h(z^k) > h^*$.  Notice from \eqref{subpr_overh} and Algorithm~\ref{algo_implement:subsolver} that $(z^{\ell+1}, z^\ell) \in \dom\,\bh$, which along with \eqref{subpr_overh} implies that $\bh(z^{\ell+1}, z^\ell) = h(z^{\ell+1})$ for all $0 \leq \ell < k$. Also, by a similar argument as used in the proof of Lemma~\ref{lem:equiv_Phi}, one has $h^* = \bh^*$, where $\bh^*$ is defined in \eqref{subpr_overh}. These, together with $h(z^0) - h^* \leq \eta$, $h(z^k) > h^*$, and the monotonicity of $\{h(z^\ell)\}$, lead to
\[
0<\bh(z^{\ell+1}, z^\ell) - \bh^* = h(z^{\ell+1})-h^* \leq h(z^0) - h^* \leq \eta \qquad \forall 0 \leq \ell < k.
\]
It then follows from \eqref{KL_subpr} that 
\beq \label{KL-z}
C(\bh(z^{\ell+1}, z^\ell) - \bh^*)^{\theta} \leq \dist(0, \partial \bh(z^{\ell+1}, z^\ell)) \qquad \forall 0 \leq \ell < k.
\eeq
In addition, notice from Theorem~\ref{thm:sub_well}(ii) and Lemma~\ref{lambda-bound} that $L_\ell \leq \bar{L}$ and $\|\lambda^\ell\|_\infty \leq A$ for all $0 \leq \ell < k$, where $\bar{L}$ and $A$ are defined in \eqref{Lk-bound} and \eqref{prop1_explicitBound}. Using these, \eqref{lem9_omega}, and  Lemma~\ref{lem:distBound}, we obtain that $\dist(0, \partial \bh(z^{\ell+1}, z^{\ell})) \leq \omega \|z^{\ell+1} - z^{\ell}\|$ for all $0 \leq \ell < k$. Also, notice from Algorithm~\ref{algo_implement:subsolver} that 
\[
r_{\ell} - r_{\ell+1} = h(z^{\ell})-h(z^{\ell+1}) \geq \frac{\beta}{2}\|z^{\ell+1}-z^{\ell}\|^2 \qquad \forall 0 \leq \ell < k.
\]
In view of these, \eqref{lem4:alpha_Cp}, and \eqref{KL-z}, one has that for all $0 \leq \ell < k$,
\begin{align}
r_{\ell} - r_{\ell+1} &\geq \frac{\beta}{2}\|z^{\ell+1}-z^{\ell}\|^2 \geq \frac{\beta}{2\omega^2} \dist^2(0, \partial \bh(z^{\ell+1}, z^{\ell})) \overset{\eqref{KL-z}}{\geq} \frac{\beta C^2}{2\omega^2} \big(\bh(z^{\ell+1}, z^{\ell}) - \bh^*\big)^{2\theta} \nn \\
&= \frac{\beta C^2}{2\omega^2} (h(z^{\ell+1}) - h^*)^{2\theta} \overset{\eqref{lem4:alpha_Cp}}{=} \alpha\,(h(z^{\ell+1}) - h^*)^{2\theta} = \alpha\,r_{\ell+1}^{2\theta}. \label{lem4_rkdiff}
\end{align}

(i) Suppose $\theta=1/2$. It then follows from \eqref{lem4_rkdiff} that $r_{\ell+1} \leq (1+\alpha)^{-1} r_{\ell}$ for all $0 \leq \ell < k$, which together with $r_0 \leq \eta$ implies that $r_k \leq r_0 (1+\alpha)^{-k} \leq \eta (1+\alpha)^{-k}$, and hence \eqref{lem4:case1} holds.

(ii) Suppose $\theta \in (1/2, 1)$. Notice from the above that $r_k > 0$, which together with the monotonicity of $\{r_\ell\}$ implies that $r_\ell>0$ for all $0 \leq \ell \leq k$. Letting $\psi(t) = \frac{1}{2\theta-1} t^{1-2\theta}$  and using the monotonicity of $\{r_\ell\}$, we have
\beq \label{r-relation}
\psi(r_{\ell+1}) - \psi(r_\ell) = \int_{r_{\ell}}^{r_{\ell+1}} \psi'(t) dt = \int_{r_{\ell+1}}^{r_\ell} t^{-2\theta} dt \geq r_\ell^{-2\theta} (r_\ell - r_{\ell+1}) \qquad \forall 0 \leq \ell < k.
\eeq
For each $0 \leq \ell < k$, we consider two separate cases below.

Case a): $r_{\ell+1}^{-2\theta} \leq 2 r_\ell^{-2\theta}$. It along with \eqref{lem4_rkdiff} and \eqref{r-relation} implies that
\[
\psi(r_{\ell+1}) - \psi(r_\ell) \geq \frac{1}{2} r_{\ell+1}^{-2\theta} (r_\ell - r_{\ell+1}) \geq \frac{1}{2} \alpha.
\]

Case b): $r_{\ell+1}^{-2\theta} > 2 r_\ell^{-2\theta}$. It leads to $r_{\ell+1}^{1-2\theta} >2^{\frac{2\theta-1}{2\theta}} r_\ell^{1 - 2\theta}$. By this, $\theta \in (1/2, 1)$, $r_\ell \leq \eta$, and the expression of $\psi$, one has that
\begin{align*}
\psi(r_{\ell+1}) - \psi(r_\ell) = \frac{1}{2\theta-1} (r_{\ell+1}^{1-2\theta} - r_\ell^{1-2\theta}) > \frac{1}{2\theta-1} \left(2^{\frac{2\theta-1}{2\theta}} - 1\right) r_\ell^{1-2\theta} \geq \frac{1}{2\theta-1} \left(2^{\frac{2\theta-1}{2\theta}} - 1\right) \eta^{1-2\theta}.
\end{align*}

Combining the above two cases and using the definition of $C'$ in \eqref{lem4:alpha_Cp}, we obtain that $\psi(r_{\ell+1}) - \psi(r_\ell) \geq C'$ for all $0 \leq \ell < k$. It follows that $\psi(r_{k}) \geq \psi(r_0)+ kC' \geq k C'$. This together with the expression of $\psi$ yields
\[
r_{k} \leq \left( \frac{1}{C' (2\theta-1)} \right)^{\frac{1}{2\theta-1}} k^{-\frac{1}{2\theta-1}}.
\]
Relation \eqref{lem4:case2} then follows from this and $r_k=h(z^k) - h^*$.
\end{proof}

We are now ready to prove Theorem \ref{thm:sub_iters}.

\begin{proof}[\textbf{Proof of Theorem~\ref{thm:sub_iters}}] 
Suppose for contradiction that Algorithm~\ref{algo_implement:subsolver} runs for more than $\overline{K}_\theta$ outer iterations. It along with \eqref{subsolver_terminate} implies that there exists some $\ell \geq \overline{K}_\theta - 1$ such that
\beq \label{thm3_contra}
\|\nabla g(z^{\ell+1}) - \nabla g(z^\ell) - L_\ell(z^{\ell+1}-z^\ell)\|^2 + 4 \Big(\sum_{j=1}^m \lambda^{\ell}_j L_{c_j} \Big)^2 \|z^{\ell+1}-z^\ell\|^2 > \tau^2.
\eeq

For notational convenience, let $r_\ell = h(z^\ell) - h^*$ and $r_{\ell+1} = h(z^{\ell+1}) - h^*$. We now show that $r_\ell \leq \beta \tau^2 /(2\omega^2)$ by considering two separate cases below.

Case a): $\theta = 1/2$. By this, \eqref{thm3:overK}, and $\ell \geq \overline{K}_\theta - 1$, one has $\ell \geq \log_{1+\alpha}\big(\frac{2\omega^2\eta}{\beta \tau^2}\big)$. It then follows from \eqref{lem4:case1} that $r_\ell \leq \eta(1+\alpha)^{-\ell} \leq \beta \tau^2/(2\omega^2)$.

Case b): $\theta \in (1/2, 1)$. Using this, \eqref{thm3:overK}, and $\ell \geq \overline{K}_\theta - 1$, we have $\ell \geq \frac{1}{C' (2\theta-1)} \big(\frac{2\omega^2}{\beta\tau^2}\big)^{2\theta-1}$. It then follows from \eqref{lem4:case2} that
\[
r_\ell \leq \Big( \frac{1}{C'(2\theta-1)} \Big)^{\frac{1}{2\theta-1}} \ell^{-\frac{1}{2\theta-1}} \leq \frac{\beta \tau^2}{2\omega^2}. 
\]
Combining these two cases, we conclude that $r_\ell \leq \beta \tau^2 / (2\omega^2)$. In addition, notice from Algorithm~\ref{algo_implement:subsolver} that 
\[
r_\ell - r_{\ell+1} = h(z^\ell)-h(z^{\ell+1}) \geq \beta\|z^{\ell+1} - z^\ell\|^2 / 2.
\]
By these relations and $r_{\ell+1} \geq 0$, one has
\beq \label{z_diff}
\|z^{\ell+1}-z^\ell\| \leq \sqrt{\frac{2(r_\ell - r_{\ell+1})}{\beta}} \leq (2/\beta)^{\frac{1}{2}} r_{\ell}^{\frac{1}{2}} \leq \tau/\omega.
\eeq
Also, using Theorem \ref{thm:sub_well}(ii), we have $0<L_{\ell} \leq \bar{L}$, where $\bar{L}$ is defined in \eqref{Lk-bound}.  This together with \eqref{lem3_b1} yields
\beq \label{grad_diff}
\|\nabla g(z^{\ell+1}) - \nabla g(z^\ell) - L_\ell(z^{\ell+1}-z^\ell)\|^2
\overset{\eqref{lem3_b1}}{\leq}  (L+L_\ell)^2  \|z^{\ell+1} - z^\ell\|^2 \leq (L+\bar{L})^2  \|z^{\ell+1} - z^\ell\|^2.
\eeq
In addition, notice from Lemma \ref{lambda-bound} that $\|\lambda^\ell\|_\infty \leq A$, where $A$ is defined in \eqref{prop1_explicitBound}. It follows from this and $\lambda^\ell\in\bR^m_+$ that 
\beq \label{aux_ineq}
\big(\sum_{j=1}^m \lambda^{\ell}_j L_{c_j} \big)^2 \leq \|\lambda^\ell\|_\infty^2 \big(\sum_{j=1}^m L_{c_j}\big)^2 \leq A^2 \big(\sum_{j=1}^m L_{c_j}\big)^2.
\eeq
Using this, \eqref{z_diff}, \eqref{grad_diff}, and the definition of $\omega$ in \eqref{lem9_omega}, we obtain that
\[
\begin{aligned}
&\|\nabla g(z^{\ell+1}) - \nabla g(z^\ell) - L_\ell(z^{\ell+1}-z^\ell)\|^2 + 4 \big(\sum_{j=1}^m \lambda^{\ell}_j L_{c_j} \big)^2 \|z^{\ell+1}-z^\ell\|^2 \\
&\overset{\eqref{grad_diff}\eqref{aux_ineq}}\leq \Big( (L+\bar{L})^2 + 4A^2 \big(\sum_{j=1}^m L_{c_j}\big)^2 \Big) \|z^{\ell+1} - z^\ell\|^2 \overset{\eqref{lem9_omega}}{=} \omega^2 \|z^{\ell+1}-z^\ell\|^2 \overset{\eqref{z_diff}}{\leq} \tau^2,
\end{aligned}
\]
which contradicts \eqref{thm3_contra}. Hence, Algorithm~\ref{algo_implement:subsolver} terminates in at most $\overline{K}_\theta$ outer iterations.

We next show that \eqref{thm3:output} holds. If $h(z^{k+1})=h^*$, \eqref{thm3:output} clearly holds. For the remainder of the proof, suppose that $h(z^{k+1}) > h^*$. By similar arguments as above, one has that $0<L_k \leq \bar{L}$, $\|\lambda^k\|_\infty \leq A$, and $\|z^{k+1}-z^k\| \leq  \tau/\omega$. It follows from these, \eqref{lem9_omega}, and \eqref{subdiff-bound} that 
\begin{align}
& \dist(0, \partial \bh(z^{k+1}, z^k)) \overset{\eqref{subdiff-bound}}{\leq} \Big( (L+L_k)^2 + 4\|\lambda^k\|_\infty^2 \big(\sum_{j=1}^m L_{c_j}\big)^2 \Big)^{\frac{1}{2}} \|z^{k+1} - z^k\| \nn \\
& \leq  \Big( (L+\bar{L})^2 + 4A^2 \big(\sum_{j=1}^m L_{c_j}\big)^2 \Big)^{\frac{1}{2}} \|z^{k+1} - z^k\|  \overset{\eqref{lem9_omega}}{=} \omega \|z^{k+1} - z^k\| \leq \tau. \label{subdiff_bound}
\end{align}
In addition, notice that $0<h(z^{k+1}) - h^* \leq h(z^0) - h^* \leq \eta$, $\bh(z^{k+1}, z^k) = h(z^{k+1})$, and $\bh^* = h^*$.  It then follows that \eqref{KL_subpr} holds with $z = z^{k+1}$ and $w=z^k$. In view of these and \eqref{subdiff_bound}, one has
\[
h(z^{k+1}) - h^* = \bh(z^{k+1}, z^k) - \bh^* \leq \big(C^{-1} \dist(0, \partial \bh(z^{k+1}, z^k))\big)^{\frac{1}{\theta}} \leq (C^{-1}\tau)^{\frac{1}{\theta}}. 
\]
Hence, \eqref{thm3:output} holds as desired.
\end{proof}

\subsection{Proof of the main results in Section~\ref{sec:FOD}}\label{sec:proof3}

In this subsection we prove Lemma \ref{lambda-uniform-bound} and Theorems~\ref{thm:outer_bound} and \ref{thm:operations}. 

\begin{proof}[\textbf{Proof of Lemma \ref{lambda-uniform-bound}}] 
Fix any $\bk \geq 0$.  Observe that $\{\lambda^{\bk,\ell}\}$ reduces to the sequence of Lagrange multipliers generated by Algorithm~\ref{algo_implement:subsolver} for solving the subproblem $\min \{-f(x^{\bk}, y) + q(y): c(y) \leq 0\}$. Also, notice from Algorithm~\ref{algo_implement:whole} that $x^{\bk}\in\cX$, which together with the definition of $G_f$ implies that
\[
\max_{y \in \cY } \big\|\nabla_y f(x^{\bk}, y)\big\| \leq G_f.
\]
Using this relation, the definition of $A_f$, and Lemma \ref{lambda-bound} with $G$ replaced by $G_f$, we obtain that $\|\lambda^{\bk,\ell}\|_1 \leq A_f$ holds for all $\ell$.  By this and the arbitrariness of $\bk$, one can see that the conclusion holds. 
\end{proof}

We now turn to the proofs of Theorems~\ref{thm:outer_bound} and \ref{thm:operations}. 
Notice that Algorithm~\ref{algo_implement:whole} shares key similarities with \cite[Algorithm~2]{lu2025first}, which is designed to solve the unconstrained nonconvex-nonconcave problem 
$\min_x \max_y \{ f(x,y) + p(x) - q(y) \}$.  In particular, Algorithm~\ref{algo_implement:whole} applies an inexact proximal gradient (IPG) method to solve $\min_x \{ F^*(x) + p(x)\}$, while \cite[Algorithm~2]{lu2025first} applies an IPG method to solve $\min_x \{ \tilde{F}^*(x) + p(x) \}$, where  
$\tilde{F}^*(x) = \max_y \{ f(x,y) - q(y)\}$. The two algorithms follow almost identical steps, differing only in how they approximate {$\nabla^\rC_\cX F^* (x^k)$ and $\nabla^\rC_\cX \tilde{F}^*(x^k)$} at a given iterate $x^k$. Specifically, Algorithm~\ref{algo_implement:whole} computes an approximation to {$\nabla^\rC_\cX F^*(x^k)$} by calling the SCP method (Algorithm~\ref{algo_implement:subsolver}) for the subproblem $\min_y \{ -f(x^k,y) + q(y) : c(y) \leq 0\}$,  whereas \cite[Algorithm~2]{lu2025first} computes an approximation to {$\nabla^\rC_\cX \tilde{F}^*(x^k)$} using a proximal gradient method for the subproblem
$\min_y \{ -f(x^k,y) + q(y) \}$.  Thanks to these close similarities, the proofs of Theorems~\ref{thm:outer_bound} and \ref{thm:operations} largely parallel those of \cite[Theorems~4 and 5]{lu2025first}. We therefore provide only a sketch of the proofs. To this end, we first present two technical lemmas. 

\begin{lemma} \label{lem:subsolver_output}
Let $\cX^\rc_\epsilon, L_{\nabla f}, C, \theta, \gamma, \sigma, \epsilon, \{\eta_\ell\}$ be given in  \eqref{X-eps-r}, Assumption \ref{ass:Lip_Smo_KL}, and Algorithm \ref{algo_implement:whole}, respectively. Suppose that $\{(x^\ell,y^\ell)\}^{k}_{\ell=0}$ are generated by Algorithm \ref{algo_implement:whole} for some $k \geq 1$ such that $x^{\ell} \in \cX^\rc_\epsilon$ for all $0\leq \ell < k$. Then, for all $0\leq \ell \leq k$,  it holds that
\begin{align}
& F^*(x^\ell) - F(x^\ell, y^\ell) \leq \min\big\{\gamma \epsilon^\sigma/2, \eta_{\ell} \big\}, \qquad \operatorname{dist}\big(y^\ell, Y^*(x^\ell)\big) \leq \frac{1}{C(1-\theta)} \min\big\{(\gamma/2)^{1-\theta} \epsilon^{\sigma (1-\theta)}, \eta_\ell^{1/2}\big\}, \label{lem9-ineq1-1}  \\ 
& \|{\nabla^\rC_\cX F^*(x^\ell)} -\nabla_x f(x^\ell, y^\ell)\| \leq \frac{L_{\nabla f}}{C(1-\theta)} \min\big\{(\gamma/2)^{1-\theta}\epsilon^{\sigma (1-\theta)},\eta_\ell^{1/2}\big\}. \label{lem9-ineq1-2} 
\end{align}
\end{lemma}

\begin{proof}
{This lemma is parallel to \cite[Lemma~9]{lu2025first}, whose proof is based on \cite[Theorems~1 and~3]{lu2025first} and \cite[Lemma~4]{lu2025first}. Note that Theorems~\ref{thm:Phi_localHolder}, \ref{thm:sub_iters}, and \ref{thm:Phi_QG} are parallel to \cite[Theorems~1 and~3]{lu2025first} and \cite[Lemma~4]{lu2025first}, respectively. Hence, the conclusion follows from Theorems~\ref{thm:Phi_localHolder}, \ref{thm:Phi_QG}, and \ref{thm:sub_iters}, together with arguments similar to those used in the proof of \cite[Lemma~9]{lu2025first}.}
\end{proof}

\begin{lemma} \label{lem:gen_converge}
Let $\epsilon>0$ be given, $M$, $\cX^\rc_\epsilon$ be defined in \eqref{X-eps-r} and \eqref{Phi_smooth_constants}, $L_f, L_{\nabla f}, C, \theta, \gamma, \sigma, \{\delta_\ell\}, \{\eta_\ell\}, \{L_\ell\}$ be given in Assumption \ref{ass:Lip_Smo_KL} and Algorithm~\ref{algo_implement:whole}, and let 
\begin{align*}
& \Delta_k := 8 \Big[\Psi(x^0) - \Psi^* + \eta_{k+1} + \sum_{\ell=0}^{k}  \Big(1+\tfrac{L_{\nabla f}^2}{(1-\theta)^2 C^2 L_\ell} \Big)\eta_\ell + \sum_{\ell=0}^{k}\frac{\delta_\ell}{2} \Big],   \\
&\underline{K}_\epsilon := \max\{k \geq 1: \Delta_k/(kL_{\lceil k/2 \rceil}) \geq \gamma^2 \epsilon^{2\sigma} / (16L_f^2)\},  \\ 
& {\overline K}_\epsilon := \max \{k \geq 0: x^k \in \cX^\rc_\epsilon \},  \\
& \ell(k) := \mathop{\arg\min}_{\lceil k/2 \rceil \leq \ell \leq k} L_\ell \|x^{\ell+1} - x^\ell\|^2.  
\end{align*}
Let $\underline{K}_\epsilon < k \leq \overline{K}_\epsilon$ be given. Suppose that $\{(x^\ell,y^\ell)\}^{k}_{\ell=0}$ are generated by Algorithm \ref{algo_implement:whole} such that $x^\ell \in \cX^\rc_\epsilon$ for all 
$0 \leq \ell \leq k$. Then we have
\[
\operatorname{dist}\big(0, \partial \Psi(x^{\ell(k)+1})\big) \leq L_{\nabla f} \sqrt{ \frac{\Delta_k}{L_{\lceil k/2 \rceil} k} } + \sqrt{ \frac{L_{k} \Delta_k}{k} } + M \Big( \frac{\Delta_k}{L_{\lceil k/2 \rceil}k} \Big)^{\frac{\nu}{2}} + (1-\theta)^{-1} C^{-1} L_{\nabla f} \eta_{\lceil k/2 \rceil}^{\frac{1}{2}}.
\]
\end{lemma}

\begin{proof}
{This lemma is parallel to \cite[Lemma~10]{lu2025first}, whose proof is based on \cite[Eqs.~(25), (72), and~(73)]{lu2025first}. Note that \eqref{Fstar-bound}, \eqref{lem9-ineq1-1}, and \eqref{lem9-ineq1-2} are parallel to \cite[Eqs.~(25), (72), and~(73)]{lu2025first}, respectively. Hence, the conclusion follows from \eqref{Fstar-bound}, \eqref{lem9-ineq1-1}, and \eqref{lem9-ineq1-2}, together with arguments similar to those used in the proof of \cite[Lemma~10]{lu2025first}.}
\end{proof}

We are now ready to provide a sketch of the proofs of Theorems~\ref{thm:outer_bound} and \ref{thm:operations}.
 
\begin{proof}[\textbf{Proof of Theorem~\ref{thm:outer_bound}}] 
{Theorem~\ref{thm:outer_bound} is parallel to \cite[Theorem~4]{lu2025first}, whose proof is based on \cite[Lemmas~9 and~10]{lu2025first}. Note that \cite[Lemmas~9 and~10]{lu2025first} are parallel to Lemmas~\ref{lem:subsolver_output} and~\ref{lem:gen_converge}, respectively. Hence, the conclusion follows from Lemmas~\ref{lem:subsolver_output} and~\ref{lem:gen_converge}, together with arguments similar to those used in the proof of \cite[Theorem~4]{lu2025first}.}
\end{proof}

\begin{proof}[\textbf{Proof of Theorem~\ref{thm:operations}}] 
{Theorem~\ref{thm:operations} is parallel to \cite[Theorem~5]{lu2025first}, whose proof is based on \cite[Theorems~3 and 4]{lu2025first}. Note that \cite[Theorems~3 and 4]{lu2025first} are parallel to Theorems~\ref{thm:sub_iters} and~\ref{thm:outer_bound}, respectively. Hence, the conclusion follows from Theorems~\ref{thm:sub_iters} and~\ref{thm:outer_bound}, together with arguments similar to those used in the proof of \cite[Theorem~5]{lu2025first}.}
\end{proof}

\end{document}